\documentstyle[amsthm,amsmath,amssymb,amsfonts,graphicx,twoside]{article}

\newtheorem{ex}{Example}
\newtheorem{th}{Theorem}
\newtheorem{dfn}{Definition}

\newtheorem{re}{Remark}
\newcommand{\slush}{\bar \circ}
\newcommand{\slh}{\bar *}

\def\o{\overline}

\makeatletter
\long\def\UR#1{\leavevmode\setbox\@tempboxa\hbox{#1}\@tempdima\fboxrule
     \advance\@tempdima \fboxsep \advance\@tempdima \dp\@tempboxa
    \hbox{\lower \@tempdima\hbox
   {\vbox{\hrule \@height \fboxrule
           \hbox{  \hskip\fboxsep
           \vbox{\vskip\fboxsep \box\@tempboxa\vskip\fboxsep}\hskip
                  \fboxsep\vrule \@width \fboxrule}%
                   }}}}
\makeatother

\makeatletter
\long\def\LR#1{\leavevmode\setbox\@tempboxa\hbox{#1}\@tempdima\fboxrule
     \advance\@tempdima \fboxsep \advance\@tempdima \dp\@tempboxa
    \hbox{\lower \@tempdima\hbox
   {\vbox{
           \hbox{  \hskip\fboxsep
           \vbox{\vskip\fboxsep \box\@tempboxa\vskip\fboxsep}\hskip
                  \fboxsep\vrule \@width \fboxrule}%
                  \hrule \@height \fboxrule}}}}
\makeatother

\makeatletter
\long\def\UL#1{\leavevmode\setbox\@tempboxa\hbox{#1}\@tempdima\fboxrule
     \advance\@tempdima \fboxsep \advance\@tempdima \dp\@tempboxa
    \hbox{\lower \@tempdima\hbox
   {\vbox{\hrule \@height \fboxrule
           \hbox{\vrule \@width \fboxrule \hskip\fboxsep
           \vbox{\vskip\fboxsep \box\@tempboxa\vskip\fboxsep}\hskip
                  \fboxsep }%
                   }}}}
\makeatother

\makeatletter
\long\def\LL#1{\leavevmode\setbox\@tempboxa\hbox{#1}\@tempdima\fboxrule
     \advance\@tempdima \fboxsep \advance\@tempdima \dp\@tempboxa
    \hbox{\lower \@tempdima\hbox
   {\vbox{
           \hbox{\vrule \@width \fboxrule \hskip\fboxsep
           \vbox{\vskip\fboxsep \box\@tempboxa\vskip\fboxsep}\hskip
                  \fboxsep }%
                  \hrule \@height \fboxrule}}}}
\makeatother

\let \ttorg \tt \def \tt{\ttorg \obeyspaces}
\newcommand{\eps}{\varepsilon}

\date{}

 \author{Louis H. Kauffman \\
 Department of Mathematics, Statistics and
 Computer Science\\
 University of Illinois at Chicago\\
 851 South Morgan St., Chicago IL 60607-7045, USA\\
 kauffman@uic.edu\\
 and\\
 Vassily O. Manturov \\
 Moscow State University\\
 Department of Mechanics and Mathematics\\
 119992, GSP-2, Leninskie Gory, MSU, Moscow, Russia\\
 vassily@manturov.mccme.ru
 }

\title{Virtual Biquandles}

\begin{document}

\maketitle

\abstract{In the present paper, we describe new approaches for
constructing virtual knot invariants. The main background of this
paper comes from formulating and bringing together the ideas of
biquandle \cite{KR}, \cite{FJK}, the virtual quandle \cite{Ma2},
the ideas of quaternion biquandles by Roger Fenn and Andrew
Bartholomew \cite{BF}, the concepts and properties of long virtual
knots \cite{Ma11}, and other ideas in the interface between
classical and virtual knot theory. In the present paper we present
a new algebraic construction of virtual knot invariants, give
various presentations of it, and study several examples. Several
conjectures and unsolved problems are presented throughout the
paper.}

\section{Introduction}

\subsection{Background}

Virtual knot theory was proposed by Kauffman
in 1996, see \cite{KaV}. The combinatorial notion of virtual knot
\footnote{In the sequel, we use the generic term ``knot'' for both
knots and links, unless otherwise specified} is defined as an
equivalence class of 4-valent plane diagrams (4-regular plane graphs with extra structure)
 where a new type of
crossing (called virtual) is allowed. This theory can be regarded
as a ``projection'' of knot theory in thickened surfaces
$S_{g}\times {\bf R}$ (as studied in \cite{JKS}). Regarded from
this point of view, virtual crossings appear as artifacts of the
diagram projection from $S_g$ to ${\bf R}^{2}$. However, the rules
for handling the virtual diagrams are motivated (in \cite{KaV}) by
the idea that one can generalize the notion of a knot diagram to
its oriented Gauss code. A Gauss code for a knot is list of
crossings encountered on traversing the knot diagram, with the
signs of the crossings indicated, and whether they are over or
under in the course of the traverse. Each crossing is encountered
twice in such a traverse, and thus the Gauss code has each
crossing label appearing twice in the list. One can define
Reidemeister moves on the Gauss codes, and thus abstract the knot
theory from its planar diagrams. Virtual knot theory is the theory
of such Gauss codes, not necessarily realizable in the plane. When
one takes such a non-realizable code, and attempts to draw a
planar diagram, virtual crossings are needed to complete the
connections in the plane. These crossings are artifacts of the
planar projection. The rules for handling virtual knot diagrams
are designed to make the representation of the virtual knot
independent of the particular choice of virtual crossings that
realizes the diagram. It turns out that these rules describe
embeddings of knots and links in thickened surfaces, stabilized by
the addition and subtraction of empty handles (i.e. the addition
and subtraction of thickened 1-handles from the surface that do
not have any part of the knot or link embedded in them)
\cite{KaV2,KaV4,Ma1,Ma8,Ma10,CKS,KUP}. \bigbreak

  Another approach to Gauss codes for knots and links is the use of Gauss diagrams as in
\cite{GPV}). In this paper by Goussarov, Polyak and Viro, the
virtual knot theory, taken as all Gauss diagrams up to
Reidemeister moves, was used to analyze the structure of Vassiliev
invariants for classical and virtual knots. In both \cite{KaV} and
\cite{GPV} it is proved that if two classical knots are equivalent
in the virtual category \cite{KUP}, then they are equivalent in
the classical category. Thus classical knot theory is properly
embedded in virtual knot theory. \bigbreak

To date, many invariants of classical knots have been generalized
for the virtual case, see \cite{GPV,KaV,KR,Ma1,Ma2,Ma8,Ma10,Saw,
SW}. In many cases, a classical invariant extends to an invariant
of virtuals. In some cases one has an invariant of virtuals that
is an extension of ideas from classical knot theory, that vanishes
or is otherwise trivial for classical knots. An example of this is
the polynomial invariant studied by Sawollek \cite{Saw}, Silver
and Williams \cite{SW} and by Kauffman and Radford
\cite{KR}\footnote{Later, we will describe another approach that
leads to the same results, \cite{Ma3}}. This invariant is produced
by the methods that give the classical Alexander polynomial, but
it is an example of a zeroth order Alexander polynomial, and is
trivial in the classical case and non-trivial in the virtual case.
Such invariants are valuable for the study of virtual knots, since
they promise the possibility of distinguishing classical from
virtual knots in key cases. Other examples of this phenomenon can
be found in \cite{Ma2,Ma5}. On the other hand, some invariants
evaluated on classical knots coincide with well known classical
knot invariants (see \cite{KaV, KaV2, KaV4, Ma3} on
generalizations of the Jones polynomial, fundamental group,
quandle and quantum link invariants). These invariants exhibit
interesting phenomena on virtual knots and links: for instance,
there exists a virtual knot $K$ with ``fundamental group''
isomorphic to ${\bf Z}$ and Jones polynomial not equal to $1$.
This phenomenon immediately implies that the knot $K$ is not
classical, and underlines the difficulty of extracting the Jones
polynomial from the fundamental group in the classical case. We
know in principle that the fundamental group, plus peripheral
information, determines the knot itself in the classical case. It
is not known how to extract the Jones polynomial from this
algebraic information. Note that, for classical knots, the quandle
is a generalzation of the
 fundamental group with a geometric interpretation. In this paper we
 consider quandles of virtual knots, defined formally in terms of their
 diagrams. However, the formally defined fundamentaly group of a virtual
 knot can be interpreted as the fundamental group of the complement of the
 virtual knot in the one-point suspension of a thickened surface where this
 knot is presented.

 Another phenomenon that does not appear in the classical case
are long knots \cite{Ma11}: if we break a virtual knot diagram at
two different points and take them to the infinity, we may obtain
two different long knots. We will discuss this subject later in
the text. \bigbreak

Beyond the fundamental group and the quandle, there are two algebraic constructions defining virtual knot
invariants: the biquandle \cite{KaV2, KR, FJK} and the virtual quandle
\cite{Ma2}. Each of them has a number of realizations (representations). In this paper, we bring these two
ideas together, and present a combined construction that
allows the extraction of new algebraic invariants of virtual knots.
\bigbreak

The general strategy of this paper is the following: One presents
an algebraically defined invariant of knots (virtual knots, long
knots, etc.) This invariant is defined axiomatically (the axioms
correspond to invariance under diagrammatic moves). In this form,
the invariant appears to have strength, but difficult to work
with. In order to manage the invariant, we may take a
representation of this invariant, say, $I$, into some category
(e.g. groups) with operations defined in terms of that
representation category. In this manner we can obtain for example,
the knot group and its generalizations. Or one can find a finite
object (or category) $G$ with operations satisfying the initial
axioms. Then the set of homomorphisms $I(K)\to G$ is an invariant
of the knot $K$. In particular, so is the {\bf cardinality} of the
set of such homomorphisms. Such homomorphisms can also be called
{\em colorings} because they correspond to colorings of the
diagram arcs (generators of $I(K)$) by elements of $G$. \bigbreak

In some cases, polynomial-type invariants emerge naturally from the algebra. In other cases it is useful
to use algebraic and polynomial-type invariants together to extract information.
\bigbreak

This paper does not pretend to describe all directions of virtual
knot theory. For instance, we say a little about the Jones
polynomial and do not describe its generalizations, the Khovanov
complex \cite{Ma12}, the surface-bracket polynomial \cite{KaV4},
and the $\Xi$-polynomial, \cite{Ma3, Ma6}. Also, we do not touch
the Vassiliev invariants for virtual links, see \cite{KaV, GPV,
Ma10}. All these concepts will be described in the book by the
authors \cite{VBOOK}. Also, a list of unsolved problems concerning
virtual knots can be found in \cite{FKM}.

\subsection{Basic definitions}

We begin with the definition of a virtual knot according to
\cite{KaV}

\begin{dfn}
A {\em virtual link diagram} is a $4$--valent graph on the plane
such that each crossing of it is either classical (i.e., one pair
of opposite edges is selected to make an overcrossing; the other
pair forms an undercrossing) or virtual (just marked by a circle).
\end{dfn}

\begin{dfn}
A {\em virtual knot} is an equivalence class of virtual knot
diagrams modulo generalized Reidemeister moves.
\end{dfn}

The set of virtual Reidemeister moves consist of:

\begin{enumerate}

\item classical Reidemeister moves:

\centering\includegraphics[width=300pt]{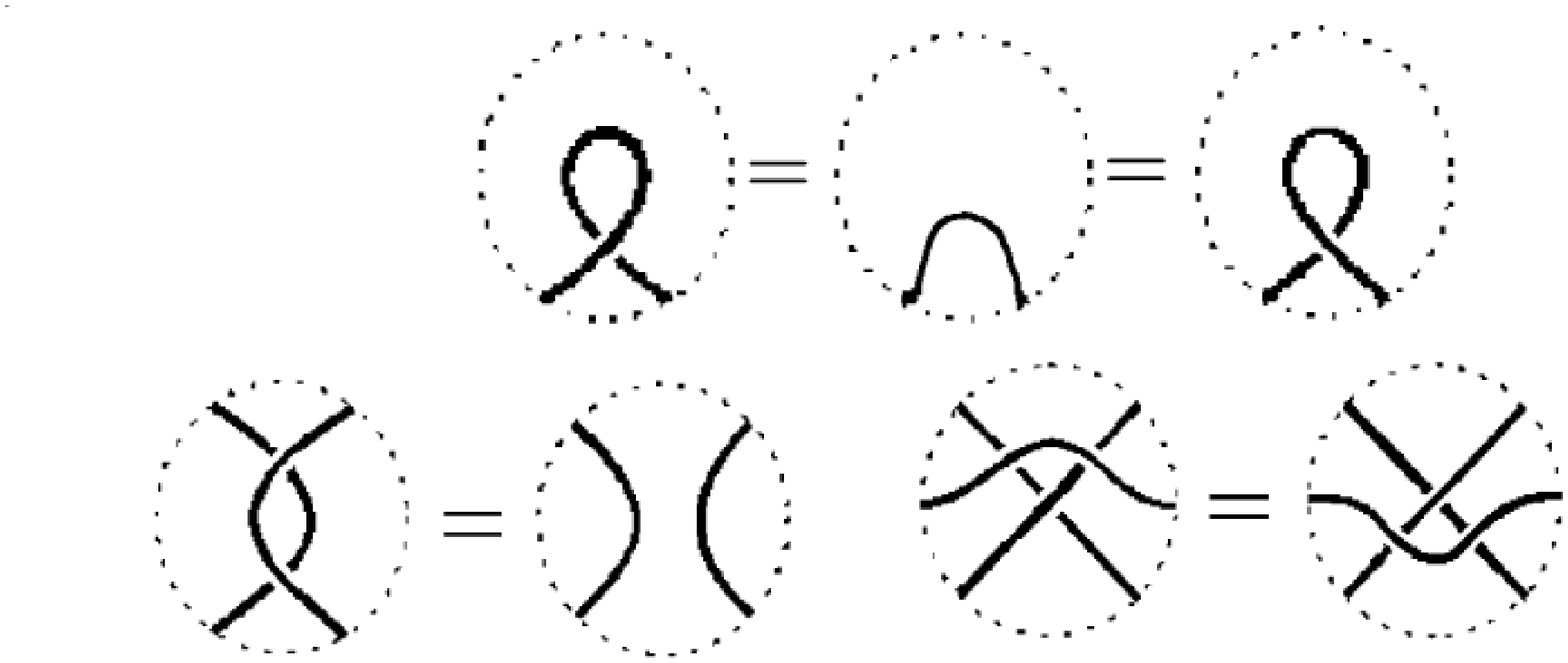}

\item virtual versions of these moves (where all classical
crossings are replaced with virtual ones)

\centering\includegraphics[width=300pt]{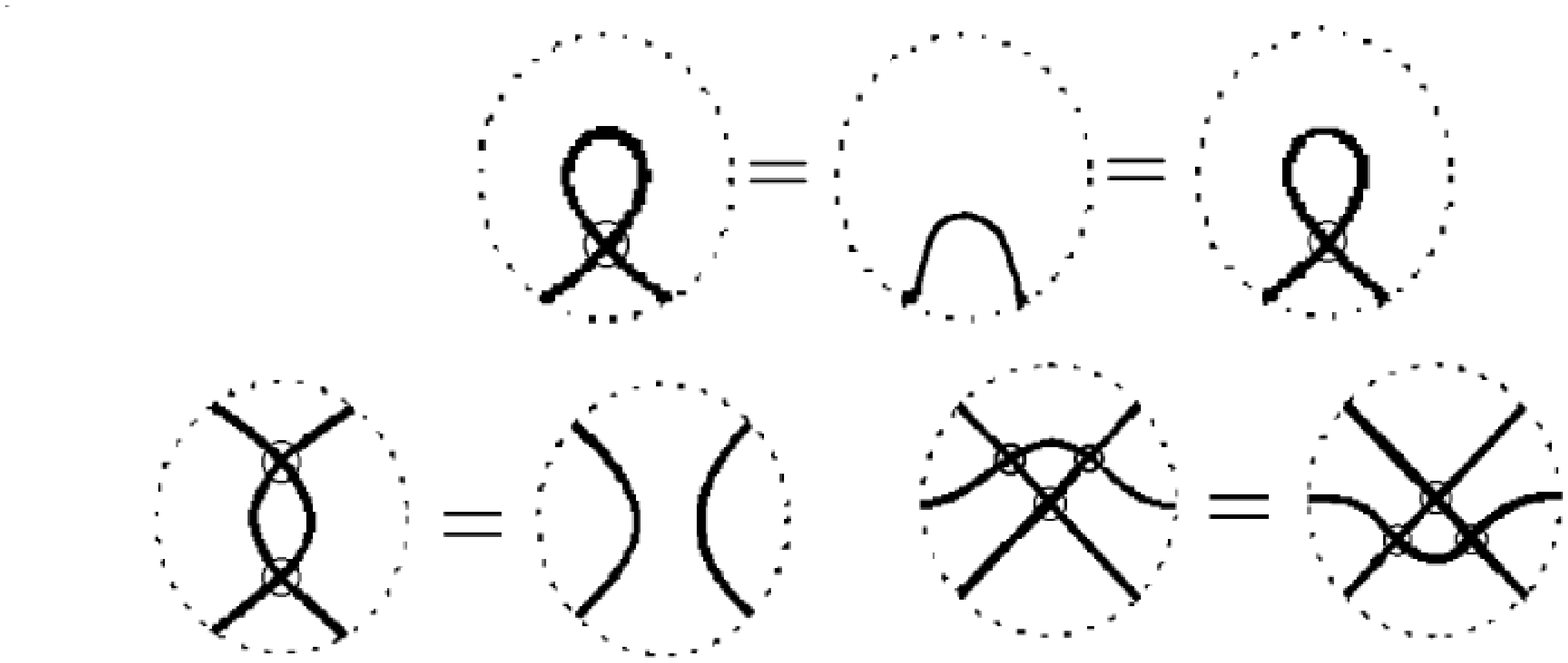}

\item and the ``semivirtual'' version of the third Reidemeister
move where two virtual crossing pass through a classical crossing.

\centering\includegraphics[width=240pt]{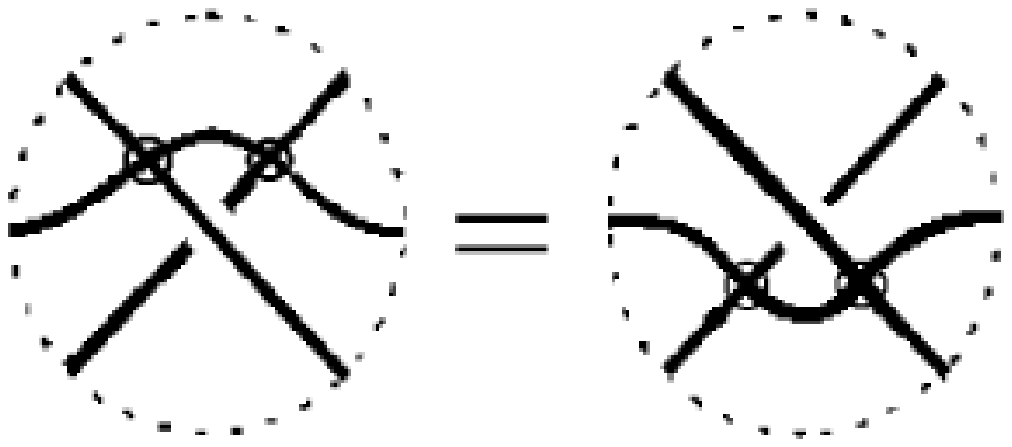}

\end{enumerate}

\begin{figure}
\centering\includegraphics[width=180pt]{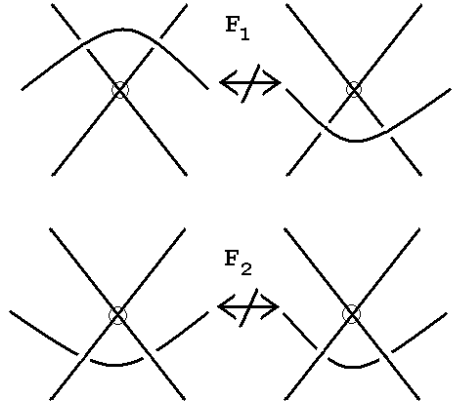}
\caption{The forbidden moves}
\label{verbot}
\end{figure}

The analogous version with two classical crossings and
one virtual crossing is forbidden. There are two versions of these forbidden
moves shown in Fig. \ref{verbot}.

It  is proved in \cite{GPV} (see also S. Nelson \cite{Nel} and T.
Kanenobu \cite{Kan}), that if we include these two moves, each
knot will be equivalent to the unknot. If we add only first of
them (F1), we will obtain what is called ``welded knots'', see
\cite{FRR} and \cite{Satoh}. In \cite{FRR} Fenn, Rimanyi and
Rourke introduced the notion of the ``welded braid group". This
version of braids is the same as virtual braids with the forbidden
move $F1$ added in braid form. Satoh in \cite{Satoh} showed that
the concept of welded knots could be interpreted in terms of
embeddings of tori in Euclidean four-space. This also leads to an
interpretation of welded braids in terms of braidings of tubes in
four-space. \bigbreak

One can easily consider {\em oriented virtual links} by endowing
diagrams with an orientation.

As with classical knots, virtual knots can be obtained as closures
of {\em virtual braids} first mentioned in the talk by Kauffman
\cite{Kau0}, and in \cite{KaV, KaV2,Ma1,KL1,KL2}. See also
Vershinin \cite{Ver} and Kamada \cite{Kam}. Virtual braids are
analogues of classical braids where some crossings are allowed to
be virtual (and marked by a circle). The corresponding $n$--strand
group has generators $\sigma_{1},\dots,\sigma_{n-1}$ (as in the
classical case) and ``virtual generators''
$\zeta_{1},\dots,\zeta_{n-1}$ with obvious relations
$\zeta_{i}^{2}=1$ (more details see in \cite{Ver, KaV, KL1}). The
closure procedure is quite analogous to that for classical knots.

\begin{re}
All virtual (and classical) knots are oriented,
unless otherwise specified.
\end{re}

The geometric approach to virtual knots is based on the following
fact: virtual knots are isotopy classes of curves in ``thickened
surfaces'' $S_{g}\times I$ up to ``stabilization'', i.e. adding
and removing handles.  See the discussion of this viewpoint given
earlier in this introduction.

\section{Basic Construction}

Here we are going to recall the construction of the quandle and
generalize it for the case of virtual knots.

\subsection{The Case of Classical Knots}

For the classical case, there is a complete algebraic invariant of
knots (complete up to reversed mirror images, and actually
complete if one adds longitudinal information to the system) first
proposed by D.Joyce \cite{Joy} and S.V.Matveev \cite{Mat}. Here is
a sketch of the quandle construction. One considers a classical
knot diagram, and encodes all arcs of it by letters. We use a
formal algebraic structure with the operation $\circ$ for which
the following condition holds. For each crossing $X$ at which the
arc $a$ (we {\bf do not look at its orientation}) lies on the
right hand of the overcrossing oriented arc $b$ and the arc $c$
lies on the left hand, we write the formal relation
\begin{equation} a\circ b=c. \label{clasrel}\end{equation} It can be checked that the invariance of
this structure (a formal algebra given by generators and
relations) under the three Reidemeister moves implies the
following algebraic axioms:

\begin{enumerate}

\item idempotence: $\forall a: a\circ a=a$

\item the existence of left inverses: $\forall b,c\; \exists ! x:
x\circ b=c$.

This $x$ is denoted by $c\slush b$.

\item right self-distributivity:
$\forall a,b,c: (a\circ b)\circ c= (a\circ c)\circ (b\circ c)$.

\end{enumerate}

An algebraic structure satisfying these axioms is called a {\em quandle}.
The quandle of an oriented knot or link is defined by generators and relations as above
(plus the imposition of these axioms).

Now, this can be easily generalized \cite{KaV} for the case of
virtual knots: we allow arcs to pass through a virtual crossing and
ignore the virtual intersections.

It is straightforward to check that this
``virtualization'' of the quandle is an invariant of virtual knots.
Since the quandle plus longitudinal information is a complete invariant of classical knots, and since
this information is preserved under virtual equivalence, one can conclude \cite{KaV,GPV}) that two classical
knot diagrams are equivalent in the virtual category if and only if they are equivalent in the classical
category.
\bigbreak

\subsection{Biquandles and a Generalized Alexander Polynomial
$G_{K}(s,t)$}

\begin{figure}
\begin{center}
{\tt    \setlength{\unitlength}{0.92pt}
\begin{picture}(280,163)
\thicklines   \put(131,162){\line(0,-1){159}}
               \put(1,1){\framebox(278,161){}}
               \put(231,104){\makebox(21,19){$a^{\o{b}} = a\, \UL{b}$}}
               \put(167,59){\makebox(20,15){$b_{\o{a}} = b\, \LL{a}$}}
               \put(233,62){\makebox(20,16){$b$}}
               \put(212,38){\makebox(18,17){$a$}}
               \put(73,106){\makebox(23,23){$a^{b} = a\, \UR{b}$}}
               \put(77,56){\makebox(32,22){$b_{a} = b\, \LR{a}$}}
               \put(10,86){\makebox(20,16){$b$}}
               \put(52,43){\makebox(18,17){$a$}}
               \put(211,89){\vector(0,1){33}}
               \put(211,43){\vector(0,1){32}}
               \put(250,82){\vector(-1,0){78}}
               \put(51,90){\vector(0,1){33}}
               \put(51,43){\vector(0,1){34}}
               \put(11,83){\vector(1,0){80}}
\end{picture}
}
\end{center}

\caption{Biquandle Relations at a Crossing}
 \label{f2}
\end{figure}
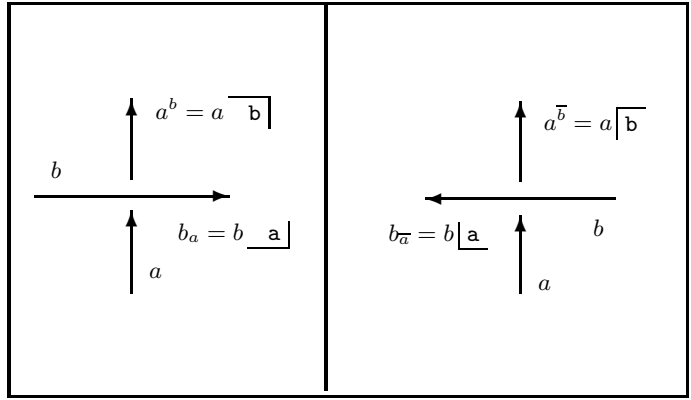

The {\em biquandle} \cite{KaV2, FJK, CS} is an algebra associated
with the diagram that is invariant (up to isomorphism) under the
generalized Reidemeister moves for virtual knots and links. The
operations in this algebra are motivated by the formation of
labels for the edges of the diagram and the intended invariance
under the moves. We will give the abstract definition of the
biquandle after a discussion of these knot theoretic issues. View
Figure \ref{f2}.  In this Figure we have shown the format for the
operations in a biquandle. The overcrossing arc has two labels,
one on each side of the crossing. In a {\em biquandle} there is an
algebra element labelling  each {\em edge} of the diagram. An edge
of the diagram corresponds to an edge of the underlying plane
graph of that diagram. \bigbreak

Let the edges oriented toward a crossing in a diagram be called
the {\em input} edges for the crossing, and the edges oriented
away from the crossing be called the {\em output} edges for the
crossing. Let $a$ and $b$ be the input edges for a positive
crossing, with $a$ the label of the undercrossing input and $b$
the label on the overcrossing input. Then in the biquandle, we
label the undercrossing output by $$c=a^{b}$$ just as in the case
of the quandle, but the overcrossing output is labeled
$$d= b_{a}.$$ We usually read $a^{b}$ as -- the undercrossing line
$a$ is acted upon by the overcrossing line $b$ to produce the
output $c=a^{b}.$ In the same way, we can read $b_{a}$ as -- the
overcossing line $b$ is operated on by the undercrossing line $a$
to produce the output $d= b_{a}.$  The biquandle labels for a
negative crossing are similar but with an overline (denoting an
operation of order two) placed on the letters. Thus in the case of
the negative crossing, we would write

  $$c=a^{\o{b}} \quad \mbox{and} \quad d=b_{\o{a}}.$$

\noindent To form the biquandle, $BQ(K)$, we take one generator
for each edge of the diagram and two relations at each crossing
(as described above). This system of generators and relations is
then regarded as encoding an algebra that is generated freely by
the biquandle operations as concatenations of these symbols and
subject to the biquandle algebra axioms.  These axioms (which we
will describe below) are a transcription in the biquandle language
of the requirement that this algebra be invariant under
Reidemeister moves on the diagram.

\noindent Another way to write this formalism for the biquandle is
as follows

$$a^{b} = a\, \UR{b}$$
$$a_{b} = a\, \LR{b}$$
$$a^{\o{b}} = a\, \UL{b}$$
$$a_{\o{b}} = a\, \LL{b}.$$

\noindent We call this the {\em operator formalism} for the
biquandle. The operator formalism has advantages when one is
performing calculations, since it is possible to maintain the
formulas on a line rather than extending them up and down the page
as in the exponential notation. On the other hand the exponential
notation has intuitive familiarity and is good for displaying
certain results. The axioms for the biquandle, are exactly the
rules needed for invariance of this structure under the
Reidemeister moves. Note that in analyzing invariance under
Reidemeister moves, we visualize representative parts of link
diagrams with biquandle labels on their edges. The primary
labeling occurs at a crossing. At a positive crossing with over
input $b$ and under input $a$, the under output is  labeled
$a\,\UR{b}$ and the over output is labeled $b\,\LR{a}\,$. At a
negative crossing with over input $b$ and under input $a$, the
under output is  labeled $a\,\UL{b}$ and the over output is
labeled $b\,\LL{a}$. At a virtual crossing there is no change in
the labeling of the lines that cross one another.

\begin{re}
Later, in this paper, we shall generalize the biquandle to include operations at the virtual crossings.
\end{re}

\begin{re} A remark is in order about the relationship of the operator
notations with the usual conventions for binary algebraic
operations.  Let $a*b = a^{b} = a\, \UR{b}.$ We are asserting that
the biquandle comes equipped with four binary operations of which
one is $a*b.$  Here is how these notations are related to the
usual parenthesizations:
\begin{enumerate}
\item $(a*b)*c = (a^{b})^{c} = a^{bc} = a\, \UR{b}\, \UR{c}$ \item
$a*(b*c) = a^{b^{c}} = a\, \UR{b\, \UR{c}}$
\end{enumerate}
\noindent From this the reader should see that the exponential and
operator notations allow us to express biquandle equation with a
minimum of parentheses.
\end{re}
\bigbreak

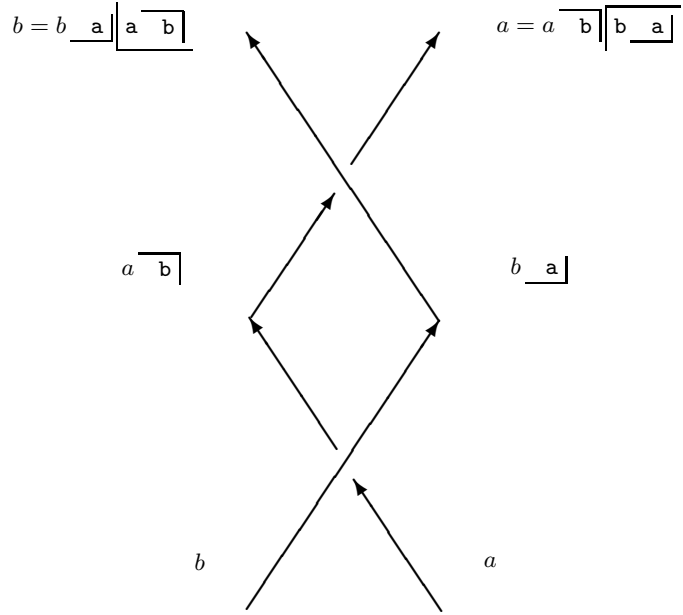
\begin{figure}
\begin{center}
{\tt    \setlength{\unitlength}{0.92pt}
\begin{picture}(323,224)
\thicklines   \put(201,222){\makebox(121,40){$a =
a\,\UR{b}\,\UL{b\,\LR{a}}$}}
               \put(41,123){\makebox(81,40){$a\,\UR{b}$}}
               \put(200,122){\makebox(82,41){$b\,\LR{a}$}}
               \put(201,2){\makebox(40,41){$a$}}
               \put(81,2){\makebox(41,41){$b$}}
               \put(1,222){\makebox(121,41){$b =
b\,\LR{a}\,\LL{a\,\UR{b}}$}}
               \put(164,186){\vector(2,3){36}}
               \put(123,123){\vector(2,3){34}}
               \put(158,69){\vector(-2,3){36}}
               \put(201,2){\vector(-2,3){36}}
               \put(200,121.5){\vector(-2,3){79}}
               \put(121,3){\vector(2,3){79}}
\end{picture}}

\end{center}

\caption{Direct Two Move } \label{f3}
\end{figure}

In Figure \ref{f3} we illustrate the effect of these conventions
and how it leads to the following algebraic transcription of the
directly oriented second Reidemeister move:

$$a = a\,\UR{b}\,\UL{b\,\LR{a}} \quad \mbox{or} \quad a = a^{b
\o{b_{a}}},$$
$$b = b\,\LR{a}\,\LL{a\,\UR{b}} \quad \mbox{or} \quad b= b_{a
\o{a^{b}}}.$$

\begin{figure}
\begin{center}
{\tt    \setlength{\unitlength}{0.92pt}
\begin{picture}(242,350)
\thicklines   \put(41,1){\makebox(160,42){$\exists x \ni x =
a\,\UR{b\,\LL{x}} \quad \mbox{,} \quad a = x\,\UL{b} \quad
\mbox{and} \quad b= b\,\LL{x}\, \LR{a} $}}
               \put(1,281){\makebox(80,42){$b\,\LL{x}\, \LR{a}$}}
               \put(26,200){\makebox(78,42){$x$}}
               \put(141,200){\makebox(80,42){$b\,\LL{x}$}}
               \put(139,80){\makebox(81,42){$x\,\UL{b}$}}
               \put(40,82){\makebox(41,39){$b$}}
               \put(141,282){\makebox(40,39){$a$}}
               \put(125,136){\vector(2,-3){34}}
               \put(83,202){\vector(2,-3){32}}
               \put(116,255){\vector(-2,-3){33}}
               \put(158,321){\vector(-2,-3){33}}
               \put(159,204){\vector(-2,3){78}}
               \put(80,84){\vector(2,3){79}}
\end{picture}}
\end{center}

\caption{Reverse Two Move } \label{f4}
\end{figure}

\noindent The reverse oriented second Reidemeister move gives a
different sort of identity, as shown in Figure \ref{f4}. For the
reverse oriented move, we must assert that given elements $a$ and
$b$ in the biquandle, then there exists an element $x$ such that
\nopagebreak $$\mbox{x} = \mbox{a}\,\UR{b\,\LL{x}} \quad \mbox{,}
\quad \mbox{a} = \mbox{x}\,\UL{b} \quad \mbox{and} \quad \mbox{b}=
\mbox{b}\,\LL{x}\, \LR{a}\,.$$ \bigbreak

By reversing the arrows in Figure \ref{f4} we obtain a second
statement for invariance under the type two move, saying the same
thing with the operations reversed: Given elements $a$ and $b$ in
the biquandle, then there exists an element $x$ such that

$$\mbox{x} = \mbox{a}\,\UL{b\,\UR{x}} \quad \mbox{,} \quad \mbox{a} = \mbox{x}\,\UR{b} \quad
\mbox{and} \quad \mbox{b}= \mbox{b}\,\UR{x}\, \LL{a}\,.$$

\noindent There is no neccessary relation between the $x$ in the
first statement and the $x$ in the second statement. \bigbreak

\begin{figure}
\begin{center}
 {\tt    \setlength{\unitlength}{0.92pt}

\begin{picture}(323,246)
\thicklines   \put(82,3){\makebox(122,40){$\exists x \ni x=a
\,\LR{x} \quad \mbox{and} \quad a = x \,\UR{a}$}}
               \put(224,162){\makebox(78,40){$x \,\UR{a}$}}
               \put(183,82){\makebox(79,40){$a \,\LR{x}$}}
               \put(61,81){\makebox(43,42){x}}
               \put(21,162){\makebox(41,41){a}}
               \put(171,148){\vector(4,3){71}}
               \put(83,82){\vector(4,3){74}}
               \put(122,43){\vector(-1,1){39}}
               \put(201,123){\vector(-1,-1){79}}
               \put(43,203){\vector(2,-1){158}}
\end{picture}}
\end{center}
\caption{First Move} \label{f5}

\end{figure}
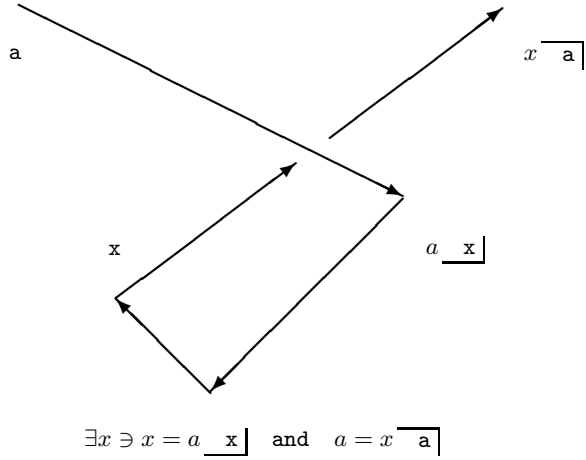

\noindent These assertions about the existence of $x$ can be
viewed as asserting the existence of fixed points for a certain
operators. In this case such an operator is $F(x)
=a\,\UR{b\,\LL{x}}\,$. It is characteristic of certain axioms in
the biquandle that they demand the existence of such fixed points.
Another example is the axiom corresponding to the first
Reidemeister move (one of them) as illustrated in Figure \ref{f5}.
This axiom states that given an element $a$ in the biquandle, then
there exists an $x$ in the biquandle such that $x=a \,\LR{x}$ and
that $a = x \,\UR{a}$. In this case the operator is $G(x) =
a\,\LR{x}\,$.

\bigbreak

It is unusual that an algebra would have axioms asserting the
existence of fixed points with respect to operations involving its
own elements. We plan to take up the study of this aspect of
biquandles in a separate publication. For now it is worth
remarking that a slight change in the axiomatic structure allows
an easy definition of the free biquandle. The idea is this:
Suppose that one has an axiom that states the existence of an $x$
such that $x=a \,\LR{x}$ for each $a.$ Then we change the
statement of the axiom by adding a new operation (unary in this
case) to the algebra, call it $Fix(a)$ such that $Fix(a) =a
\,\LR{Fix(a)}.$ Existence of the fixed point follows from this
property of the new operation, and we can describe the free
biquandle on a set by taking all finite biquandle expressions in
the elememts of the set, modulo these revised axioms for the
biquandle. \bigbreak

The biquandle relations for invariance under the third
Reidemeister move are shown in Figure \ref{f6}. The version of the
third Reidemeister move shown in this figure yields the algebraic
relations:

$$a\, \UR{b} \,\UR{c} = a\, \UR{c\, \LR{b}} \,\UR{b\, \UR{c}} \quad
\mbox{or} \quad a^{b c} = a^{c_{b} b^{c}},$$

$$c\, \LR{b} \,\LR{a} = c\, \LR{a\, \UR{b}} \,\LR{b\, \LR{a}} \quad
\mbox{or} \quad c_{b a} = c_{a^{b} b_{a}},$$

$$b\, \LR{a} \, \UR{c\, \LR{a\, \UR{b}}} = b\, \UR{c}\, \LR{a\,
\UR{c\, \LR{b}}} \quad \mbox{or} \quad
  (b_{a})^{c_{a^{b}}} = (b^{c})_{a^{c_{b}}}.$$

\begin{figure}
\begin{center}
{\tt    \setlength{\unitlength}{0.92pt}
\begin{picture}(173,371)
\thicklines   \put(120,63){\vector(1,1){26}}
               \put(87,29){\vector(1,1){25}}
               \put(63,53){\vector(1,-1){24}}
               \put(26,92){\vector(1,-1){28}}
               \put(121,314){\vector(1,-1){25}}
               \put(84,350){\vector(1,-1){28}}
               \put(59,324){\vector(1,1){25}}
               \put(26,290){\vector(1,1){25}}
               \put(147,148){\vector(-1,-1){120}}
               \put(27,148){\vector(1,-1){56}}
               \put(92,82){\vector(1,-1){55}}
               \put(91,283){\vector(1,-1){55}}
               \put(26,349){\vector(1,-1){56}}
               \put(146,349){\vector(-1,-1){120}}
\thinlines    \put(196,2){\makebox(23,23){$b\, \UR{c}\, \LR{a\,
\UR{c\, \LR{b}}}$}}
               \put(76,5){\makebox(22,21){$a\, \UR{c\, \LR{b}}$}}
               \put(0,1){\makebox(25,24){$c\, \LR{b} \,\LR{a}$}}
               \put(102,72){\makebox(26,24){$b\, \UR{c}$}}
               \put(48,72){\makebox(23,24){$c\, \LR{b}$}}
               \put(196,77){\makebox(22,24){$a\, \UR{c\, \LR{b}}
\,\UR{b\, \UR{c}}$}}
               \put(149,138){\makebox(23,25){$c$}}
               \put(5,139){\makebox(21,23){$b$}}
               \put(4,79){\makebox(20,21){$a$}}
               \put(144,200){\makebox(22,23){$b\, \LR{a} \, \UR{c\,
\LR{a\, \UR{b}}}$}}
               \put(6,203){\makebox(25,23){$c\, \LR{a\, \UR{b}}
\,\LR{b\, \LR{a}}$}}
               \put(164,279){\makebox(22,25){$a\, \UR{b} \,\UR{c}$}}
               \put(104,280){\makebox(21,23){}}
               \put(48,281){\makebox(23,21){$b\, \LR{a}$}}
               \put(149,340){\makebox(22,21){c}}
               \put(74,349){\makebox(21,21){$a\, \UR{b}$}}
               \put(3,349){\makebox(23,21){$b$}}
               \put(1,282){\makebox(20,20){$a$}}
\end{picture}}

\end{center}
\caption{Third Move}
\label{f6}
\end{figure}

\noindent The reader will note that if we replace the diagrams of
Figure \ref{f6} with diagrams with all negative crossings then we
will get a second triple of equations identical to the above
equations but with all right operator symbols replaced by the
corresponding left operator symbols (equivalently -- with all
exponent literals replaced by their barred versions).  Here are
the operator versions of these equations. We refrain from writing
the exponential versions because of the prolixity of barred
variables.

$$a\, \UL{b} \,\UL{c} = a\, \UL{c\, \LL{b}} \,\UL{b\, \UL{c}},$$

$$c\, \LL{b} \,\LL{a} = c\, \LL{a\, \UL{b}} \,\LL{b\, \LL{a}},$$

$$b\, \LL{a} \, \UL{c\, \LL{a\, \UL{b}}} = b\, \UL{c}\, \LL{a\,
\UL{c\, \LL{b}}}.$$

\vspace{3mm}

We now have a complete set of axioms, for it is a fact (see, e.g.
\cite{KNOTS,Ma1}) that the third Reidemeister move with the orientation
shown in Figure \ref{f6} and either all positive crossings (as
shown in that Figure) or all negative crossings, is sufficient to
generate all the other cases of third Reidemeister move just so
long as we have both oriented forms of the second Reidemeister
move. Consequently, we can now give the full definition of the
biquandle. \bigbreak

\noindent {\bf Definition.}  A {\em biquandle} $B$ is a set with
four binary operations indicated by the conventions we have
explained above:  $a^{b} \,\mbox{,} \, a^{\o{b}} \, \mbox{,} \,
a_{b} \,\mbox{,} \, a_{\o{b}}.$ We shall refer to the operations
with barred variables as the {\em left} operations and the
operations without barred variables as the {\em right} operations.
The biquandle is closed under these operations and the following
axioms are satisfied:\nopagebreak
\begin{enumerate}
\item   For any elements $a$ and $b$ in $B$ we have

$$a = a^{b \o{b_{a}}}  \quad \mbox{and} \quad  b= b_{a \o{a^{b}}}
\quad \mbox{and}$$

$$a = a^{\o{b}b_{\o{a}}}  \quad \mbox{and} \quad  b= b_{\o{a}
a^{\o{b}}}.$$

\item  Given elements $a$ and $b$ in $B$, then there exists an
element $x$ such that

$$x = a^{b_{\o{x}}} \mbox{,} \quad a = x^{\o{b}} \quad \mbox{and}
\quad b= b_{\o{x}a}.$$

\noindent Given elements $a$ and $b$ in $B$, then there exists an
element $x$ such that

$$x = a^{\o{b_{x}}} \mbox{,} \quad a = x^{b} \quad \mbox{and} \quad
b= b_{x\o{a}}.$$

\item For any $a$ , $b$ , $c$ in $B$ the following equations hold
and the same equations hold when all right operations are replaced
in these equations by left operations.

$$a^{b c} = a^{c_{b} b^{c}} \mbox{,} \quad c_{b a} = c_{a^{b} b_{a}}
\mbox{,} \quad (b_{a})^{c_{a^{b}}} = (b^{c})_{a^{c_{b}}}.$$

\item Given an element $a$ in $B$, then there exists an $x$ in the
biquandle such that $x=a_{x}$ and $a = x^{a}.$ Given an element
$a$ in $B$, then there exists an $x$ in the biquandle such that
$x=a^{\o{x}}$ and $a = x_{\o{a}}.$
\end{enumerate}

These axioms are transcriptions of the Reidemeister moves. The
first axiom transcribes the directly oriented second Reidemeister
move. The second axiom transcribes the reverse oriented
Reidemeister move. The third axiom transcribes the third
Reidemeister move as we have described it in Figure \ref{f6}. The
fourth axiom transcribes the first Reidemeister move. Much more
work is needed in exploring these algebras and their applications
to knot theory. \bigbreak

\subsection{The Alexander Biquandle}

\noindent In order to realize a specific example of a biquandle
structure, suppose that

$$a\,\UR{b} = ta + vb$$
$$a\,\LR{b} =sa + ub$$

\noindent where  $a$,$b$,$c$ are elements of a module $M$ over a
ring $R$ and $t$,$s$,$v$ and $u$ are in $R$. We use invariance
under the Reidemeister moves to determine relations among these
coefficients. \vspace{3mm}

Taking the equation for the third Reidemeister move discussed
above, we have

$$a \UR{b} \,\UR{c} = t(ta+vb) + vc = t^{2}a + tvb +vc$$
$$a \UR{c \LR{b}} \,\UR{b \UR{c}} = t(ta + v(sc + ub)) + v(tb + vc)$$
$$= t^{2}a + tv(u+1)b + v(ts+v)c.$$

\noindent From this we see that we have a solution to the equation
for the third Reidemeister move if $u=0$ and $v=1-st$. Assuming
that $t$ and $s$ are invertible, it is not hard to see that the
following equations not only solve this single Reideimeister move,
but they give a biquandle structure, satisfying all the axioms.

$$a\,\UR{b} = ta + (1-st)b \, \mbox{,} \quad a\,\LR{b} = sa$$
$$a\,\UL{b} = t^{-1}a + (1-s^{-1}t^{-1})b \, \mbox{,} \quad a\,\LL{b}
= s^{-1}a.$$

\noindent Thus we have a simple generalization of the Alexander
quandle and we shall refer to this structure, with the equations
given above, as the {\em Alexander Biquandle}. \vspace{3mm}

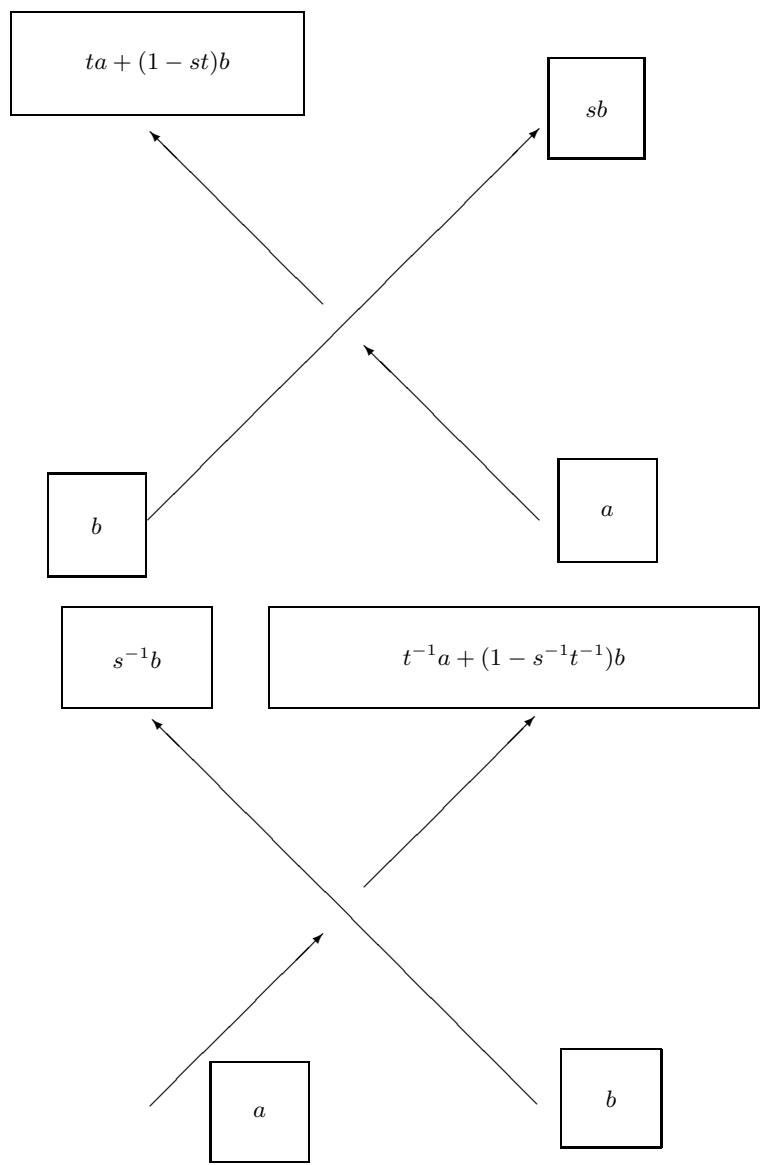
\begin{figure}
\begin{center}
{\tt    \setlength{\unitlength}{0.92pt}
\begin{picture}(309,475)
\thinlines    \put(217,25){\vector(-1,1){158}}
               \put(58,24){\vector(1,1){71}}
               \put(146,114){\vector(1,1){70}}
               \put(83,1){\framebox(40,41){$a$}}
               \put(227,7){\framebox(41,40){$b$}}
               \put(22,188){\framebox(61,41){$s^{-1}b$}}

\put(107,188){\framebox(201,41){$t^{-1}a+(1-s^{-1}t^{-1})b$}}
               \put(1,432){\framebox(120,42){$ta+(1-st)b$}}
               \put(222,414){\framebox(39,41){$sb$}}
               \put(16,242){\framebox(40,42){$b$}}
               \put(226,248){\framebox(40,42){$a$}}
               \put(129,354){\vector(-1,1){71}}
               \put(218,265){\vector(-1,1){72}}
               \put(57,265){\vector(1,1){161}}
\end{picture}}

\end{center}

\caption{Alexander Biquandle Labeling at a Crossing} \label{f7}
\end{figure}

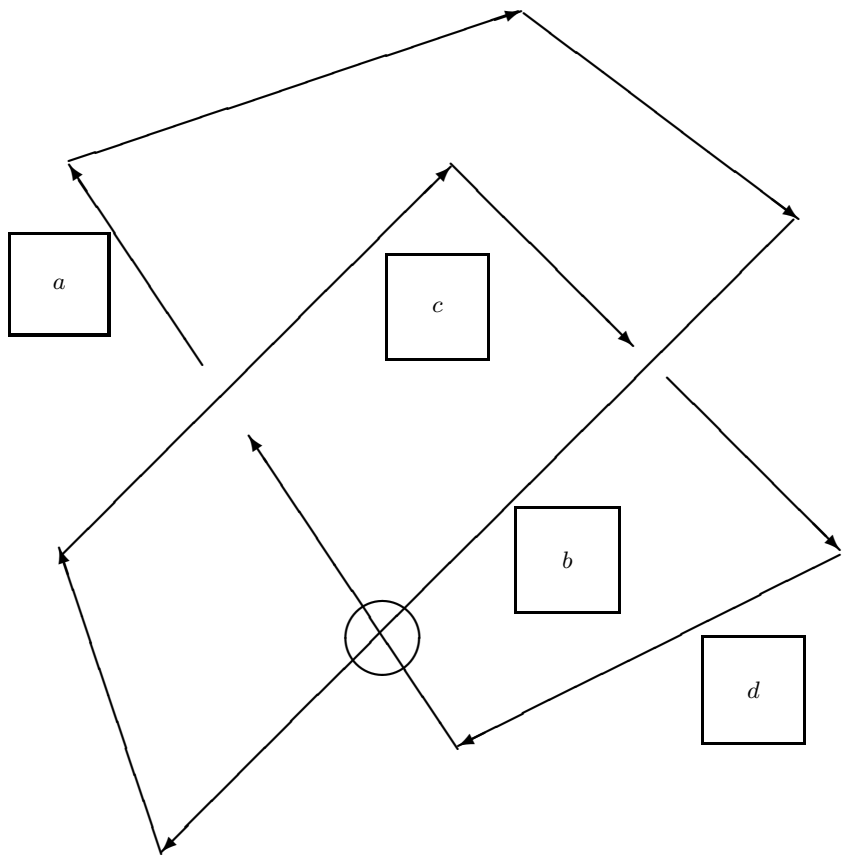
\begin{figure}
\begin{center}

{\tt    \setlength{\unitlength}{0.92pt}
\begin{picture}(344,352)
\thicklines   \put(286,48){\framebox(41,43){$d$}}
               \put(156,206){\framebox(41,42){$c$}}
               \put(209,102){\framebox(42,42){$b$}}
               \put(1,216){\framebox(40,41){$a$}}
               \put(154,91){\circle{30}}
               \put(271,198){\vector(1,-1){71}}
               \put(182,286){\vector(1,-1){75}}
               \put(80,203){\vector(-2,3){55}}
               \put(185,45){\vector(-2,3){86}}
               \put(63,2){\vector(-1,3){42}}
               \put(323,263){\vector(-1,-1){260}}
               \put(212,348){\vector(4,-3){113}}
               \put(25,287){\vector(3,1){186}}
               \put(342,125){\vector(-2,-1){157}}
               \put(22,125){\vector(1,1){160}}
\end{picture}}

\end{center}
\caption{A Virtual Knot Fully Labeled}

\label{f8}
\end{figure}

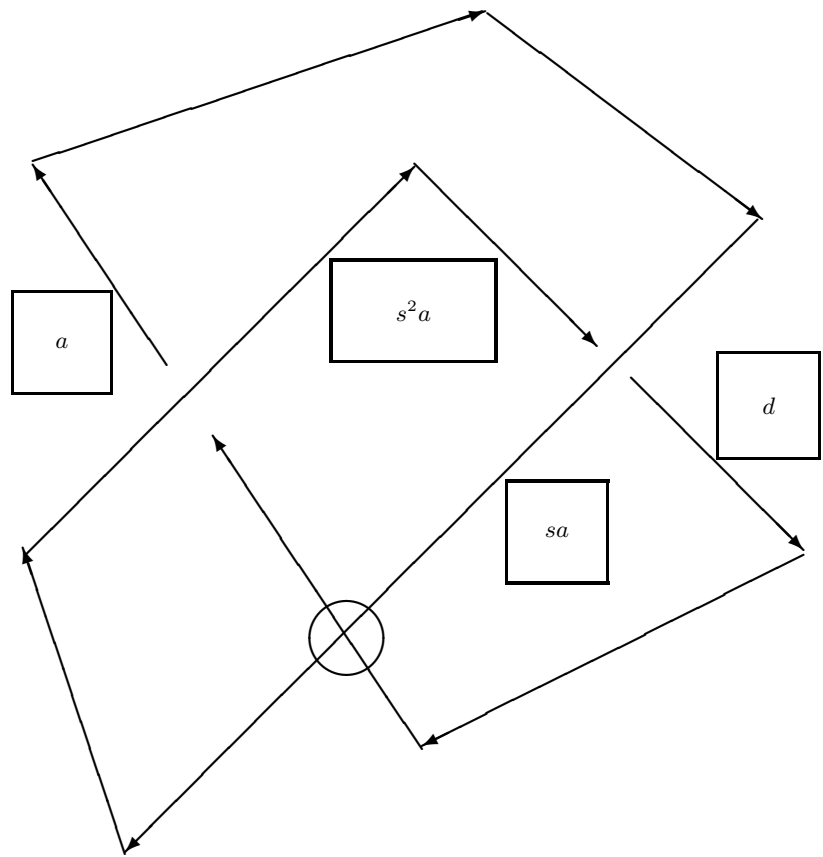
\begin{figure}
\begin{center}
 {\tt\setlength{\unitlength}{0.92pt}
\begin{picture}(333,352)
\thicklines   \put(132,205){\framebox(67,41){$s^{2}a$}}
               \put(204,114){\framebox(41,41){$sa$}}
               \put(291,165){\framebox(41,43){$d$}}
               \put(1,192){\framebox(40,41){$a$}}
               \put(138,91){\circle{30}}
               \put(255,198){\vector(1,-1){71}}
               \put(166,286){\vector(1,-1){75}}
               \put(64,203){\vector(-2,3){55}}
               \put(169,45){\vector(-2,3){86}}
               \put(47,2){\vector(-1,3){42}}
               \put(307,263){\vector(-1,-1){260}}
               \put(196,348){\vector(4,-3){113}}
               \put(9,287){\vector(3,1){186}}
               \put(326,125){\vector(-2,-1){157}}
               \put(6,125){\vector(1,1){160}}
\end{picture}}
\end{center}

\caption{A Virtual Knot with Lower Operations Labeled}
\label{f9}
\end{figure}

Just as one can define the Alexander Module of a classical knot,
we have the Alexander Biquandle of a virtual knot or link,
obtained by taking one generator for each {\em edge} of the knot
diagram and taking the module relations in the above linear form.
Let $ABQ(K)$ denote this module structure for an oriented link
$K$. That is, $ABQ(K)$ is the module generated by the edges of the
diagram, modulo the submodule generated by the relations. This
module then has a biquandle structure specified by the operations
defined above for an Alexnder Biquandle.  We first construct the
module and then note that it has a biquandle structure. See
Figures \ref{f7},\ref{f8} and \ref{f9} for an illustration of the
Alexander Biquandle labelings at a crossing. \vspace{3mm}

For example, consider the virtual knot in Figure \ref{f8}.  This
knot gives rise to a biquandle with generators $a$,$b$,$c$,$d$ and
relations

$$a= d\, \UR{b} \, \mbox{,} \quad c= b\, \LR{d} \, \mbox{,} \quad d=
c\, \UR{a} \, \mbox{,} \quad b= a\, \LR{c}.$$

\noindent writing these out in $ABQ(K)$, we have

$$a = td + (1-st)b \, \mbox{,} \quad c=sb \, \mbox{,} \quad
d=tc+(1-st)a \, \mbox{,} \quad b=sa.$$

\noindent eliminating $c$ and $b$ and rewriting, we find

$$a=td+(1-st)sa$$
$$d=ts^{2}a+(1-st)a$$

Note that these relations can be written directly from the diagram
as indicated in Figure \ref{f9} if we perform the lower biquandle
operations directly on the diagram.  This is the most convenient
algorithm for producing the relations. \smallbreak

\noindent We can write these as a list of relations

$$(s-s^{2}t-1)a +td = 0$$
$$(s^{2}t+1-st)a -d = 0$$

\noindent for the Alexander Biquandle as a module over
$Z[s,s^{-1}, t, t^{-1}].$ The relations can be expressed concisely
with the matrix of coefficients of this system of equations:

  $$M =  \left[
\begin{array}{cc}
      s-s^{2}t-1 & t  \\
      (s^{2}t+1-st) & -1
\end{array}
\right]. $$

\noindent The determinant of $M$ is, up to multiples of $\pm
s^{i}t^{j}$ for integers $i$ and $j$, an invariant of the virtual
knot or link $K$. We shall denote this determinant by $G_{K}(s,t)$
and call it the {\em generalized Alexander polynomial for $K$}. A
key fact about $G_{K}(s,t)$ is that {\em $G_{K}(s,t)=0$ if $K$ is
equivalent to a classical diagram}. This is seen by noting that in
a classical diagram one of the relations will be a consequence of
the others. \smallbreak

\noindent In this case we have

$$G_{K} = (1-s) +(s^{2}-1)t + (s-s^{2})t^{2},$$
\noindent which shows that the knot in question is non-trivial and
non-classical. \bigbreak

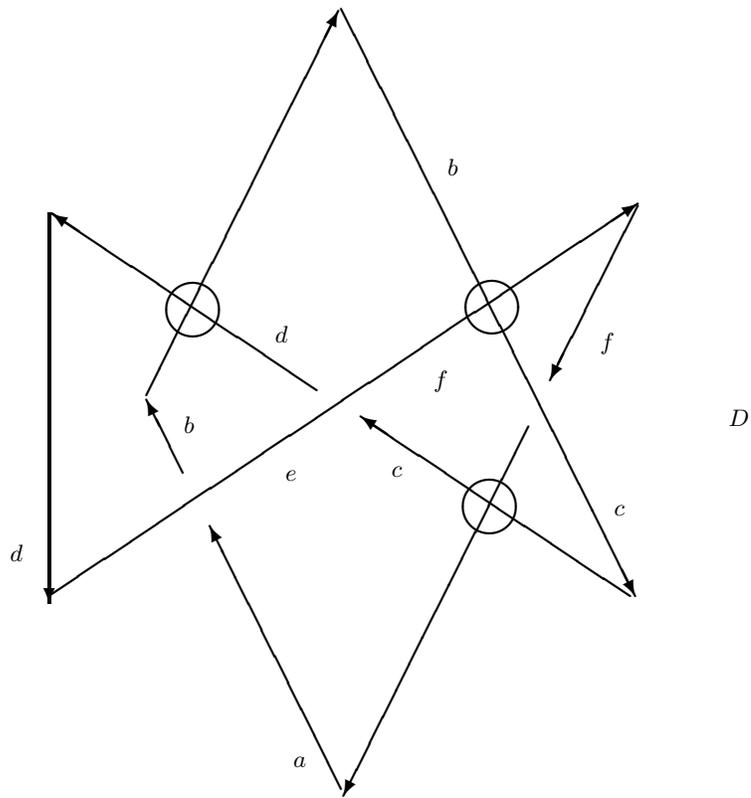
\begin{figure}
\begin{center}

{\tt    \setlength{\unitlength}{0.92pt}
\begin{picture}(326,328)
\thinlines    \put(292,137){\makebox(33,41){$D$}} \thicklines
\put(25,242){\vector(0,-1){161}}
               \put(26,85){\vector(3,2){241}}
               \put(65,167){\vector(1,2){79}}
               \put(145,326){\vector(1,-2){121}}
               \put(145,5){\vector(-1,2){54}}
               \put(80,135){\vector(-1,2){15}}
               \put(264,84){\vector(-3,2){111}}
               \put(135,169){\vector(-3,2){109}}
               \put(267,245){\vector(-1,-2){36}}
               \put(222,154){\vector(-1,-2){76}}
               \put(372,291){\circle{0}}
               \put(84,202){\circle{22}}
               \put(207,203){\circle{22}}
               \put(206,121){\circle{22}}
               \put(117,4){\makebox(22,24){$a$}}
               \put(73,144){\makebox(19,21){$b$}}
               \put(180,248){\makebox(22,25){$b$}}
               \put(249,106){\makebox(21,27){$c$}}
               \put(160,123){\makebox(16,25){$c$}}
               \put(111,181){\makebox(19,22){$d$}}
               \put(1,90){\makebox(21,24){$d$}}
               \put(116,124){\makebox(17,20){$e$}}
               \put(179,164){\makebox(14,17){$f$}}
               \put(244,177){\makebox(21,22){$f$}}
\end{picture}}
\end{center}

\caption{Unit Jones, Integer Fundamental Group}

\label{f10}
\end{figure}

Here is another example of the use of this polynomial.  Let $D$
denote the diagram in Figure \ref{f10}. It is not hard to see that
this virtual knot has unit Jones polynomial, and that the
fundamental group is isomorphic to the integers.  The biquandle
does detect the knottedness of $D$. The relations are

$$ a\, \UR{d} = b, \, d\, \LR{a} = e \, \mbox{,} \quad c\, \UR{e} =
d, \, e\,\LR{c} = f \, \mbox{,} \quad f\, \UL{b} = a, \, b\,
\LL{f} = c$$

\noindent from which we obtain the relations (eliminating $c$, $e$
and $f$).

$$b = ta +(1-tv)d \, \mbox{,} \quad d = ts^{-1}b + (1-ts)sd \,
\mbox{,} \quad a = t^{-1}s^{2}d + (1-t^{-1}s^{-1})b \, .$$

\noindent The determinant of this system is the generalized
Alexander polynomial for $D$:

$$t^{2}(s^{2}-1) + t(s^{-1}+1-s-s^{2}) + (s-s^{2}).$$

\noindent This proves that $D$ is a non-trivial virtual knot.

In fact the polynomial that we have computed is the same as the
polynomial invariant of virtuals of Sawollek \cite{Saw} and
defined by an alternative method by Silver and Williams \cite{SW}
and, in a third way, by Manturov \cite{Ma3} and given a state sum
formulation by Kauffman and Radford \cite{KR}. Sawollek defines a
module structure essentially the same as our Alexander Biquandle.
Silver and Williams first define a group. The Alexander Biquandle
proceeds from taking the abelianization of the Silver-Williams
group. Manturov uses the virtual quandle construction we shall
describe in the next subsection.\vspace{3mm}

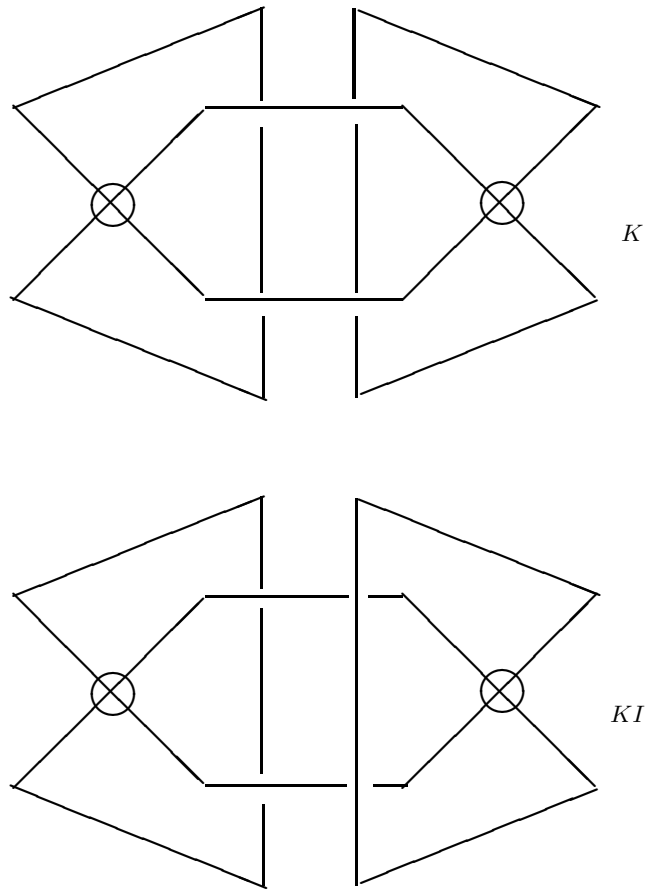
\begin{figure}
\begin{center}
{\tt    \setlength{\unitlength}{0.92pt}
\begin{picture}(272,368)
\thicklines   \put(243,64){\makebox(25,21){$KI$}}
               \put(244,260){\makebox(27,25){$K$}}
               \put(165,45){\line(-1,0){14}}
               \put(82,45){\line(1,0){58}}
               \put(163,123){\line(-1,0){14}}
               \put(82,123){\line(1,0){59}}
               \put(144,163){\line(0,-1){159}}
               \put(3,124){\line(1,-1){78}}
               \put(3,44){\line(1,1){78}}
               \put(163,124){\line(1,-1){79}}
               \put(163,44){\line(1,1){80}}
               \put(105,164){\line(0,-1){38}}
               \put(105,118){\line(0,-1){68}}
               \put(106,37){\line(0,-1){33}}
               \put(145,163){\line(5,-2){99}}
               \put(44,83){\circle{18}}
               \put(204,84){\circle{18}}
               \put(3,123){\line(5,2){103}}
               \put(2,45){\line(5,-2){105}}
               \put(146,5){\line(5,2){97}}
               \put(146,206){\line(5,2){97}}
               \put(2,246){\line(5,-2){105}}
               \put(3,324){\line(5,2){103}}
               \put(204,285){\circle{18}}
               \put(44,284){\circle{18}}
               \put(145,364){\line(5,-2){99}}
               \put(144,238){\line(0,-1){33}}
               \put(144,317){\line(0,-1){69}}
               \put(143,365){\line(0,-1){37}}
               \put(106,238){\line(0,-1){33}}
               \put(105,316){\line(0,-1){68}}
               \put(105,365){\line(0,-1){38}}
               \put(82,245){\line(1,0){81}}
               \put(82,324){\line(1,0){81}}
               \put(163,245){\line(1,1){80}}
               \put(163,325){\line(1,-1){79}}
               \put(3,245){\line(1,1){78}}
               \put(3,325){\line(1,-1){78}}
\end{picture}}
\end{center}
\caption{The Knot $K$ and the Kishino Diagram $KI$ }

\label{f11}
\end{figure}

We end this discussion of the Alexander Biquandle with two
examples that show clearly its limitations. View Figure \ref{f11}.
In this Figure we illustrate two diagrams labeled $K$ and $KI.$ It
is not hard to calculate that both $G_{K}(s,t)$ and $G_{KI}(s,t)$
are equal to zero. However, The Alexander Biquandle of $K$ is
non-trivial -- calculation shows that it is isomorphic to the free
module over $Z[s, s^{-1},t, t^{-1}]$ generated by elements $a$ and
$b$ subject to the relation $(s^{-1} - t -1)(a-b) =0.$ Thus $K$
represents a non-trivial virtual knot. This shows that it is
possible for a non-trivial virtual diagram to be a connected sum
of two trivial virtual diagrams, and it shows that the Alexander
Biquandle can sometimes be more powerful than the polynomial
invariant $G.$ However, the diagram $KI$ also has trivial
Alexander Biquandle.  This knot is proved to be knotted virtual by
Kishino \cite{KS}, Bartholomew and Fenn \cite{BF}, Manturov
\cite{Ma8}, Kadokami \cite{Kad}, Dye and Kauffman \cite{KaV4} and
others (see also the first problem in \cite{FKM}. \bigbreak

\subsection{Virtual Quandles}
                                  \begin{dfn}
A {\em virtual quandle} is a quandle $(M,\circ)$
endowed with a unary operation $f$ such that:

\begin{enumerate}

\item $f$ is invertible, the inverse operation is denoted by $f^{-1}$;

\item $\circ$ is distributive with respect to $f$:

\begin{equation} \forall a,b\in M: f(a)\circ f(b)=f(a\circ
b).\label{eqeq}\end{equation}

\end{enumerate}
\end{dfn}

\begin{re}
The equation \ref{eqeq} implies that for all $a,b\in M$:
$$f^{-1}(a)\circ f^{-1}(b)=f^{-1}(a\circ b),$$
$$f(a)\slush f(b)=f(a\slush b),$$ and
$$f^{-1}(a)\slush f^{-1}(b)=f^{-1}(a\slush b).$$
\end{re}

Given a virtual link diagram $L$, we construct its
virtual quandle $Q(L)$ as follows.

\begin{figure}
\begin{center}
\unitlength 1mm
\linethickness{0.8pt}
\begin{picture}(94.00,43.00)
\put(26.67,14.33){\vector(1,1){0.2}}
\bezier{56}(15.00,10.00)(23.33,9.67)(26.67,14.33)
\bezier{100}(26.67,14.00)(35.00,19.00)(30.00,33.33)
\bezier{72}(30.00,33.33)(23.00,39.67)(17.33,32.67)
\bezier{92}(17.33,32.67)(11.33,25.00)(13.33,11.67)
\bezier{16}(13.33,11.67)(13.33,9.67)(15.00,10.00)
\bezier{52}(22.00,12.33)(18.00,16.67)(19.33,23.67)
\bezier{60}(19.33,23.67)(21.33,31.00)(27.67,34.00)
\bezier{96}(30.67,34.67)(38.00,36.33)(39.00,20.33)
\put(24.33,10.33){\vector(-2,1){0.2}}
\bezier{148}(39.00,21.67)(39.00,2.00)(24.33,10.33)
\put(76.67,14.33){\vector(1,1){0.2}}
\bezier{56}(65.00,10.00)(73.33,9.67)(76.67,14.33)
\bezier{72}(80.00,33.33)(73.00,39.67)(67.33,32.67)
\bezier{92}(67.33,32.67)(61.33,25.00)(63.33,11.67)
\bezier{16}(63.33,11.67)(63.33,9.67)(65.00,10.00)
\bezier{52}(72.00,12.33)(68.00,16.67)(69.33,23.67)
\bezier{60}(69.33,23.67)(71.33,31.00)(77.67,34.00)
\put(74.33,10.33){\vector(-2,1){0.2}}
\bezier{148}(89.00,21.67)(89.00,2.00)(74.33,10.33)
\bezier{80}(89.00,21.67)(89.00,36.67)(83.67,35.33)
\bezier{16}(83.67,35.33)(81.67,36.00)(80.00,35.33)
\bezier{52}(76.67,14.33)(80.33,18.67)(88.00,19.00)
\bezier{32}(89.67,19.33)(94.00,20.00)(93.00,23.00)
\bezier{20}(93.00,23.00)(91.33,26.33)(89.67,26.67)
\bezier{44}(87.67,27.67)(82.33,29.33)(80.00,33.67)
\put(51.33,22.00){\makebox(0,0)[cc]{$\stackrel{\Omega_2}{\to}$}}
\put(22.67,38.00){\makebox(0,0)[cc]{$\downarrow$}}
\put(22.67,43.00){\makebox(0,0)[cc]{circular long arc}}
\put(77.33,40.33){\makebox(0,0)[cc]{normal diagram}}
\end{picture}

\end{center}

\caption{Reconstructing a link diagram in a proper way}
\label{deform}
\end{figure}
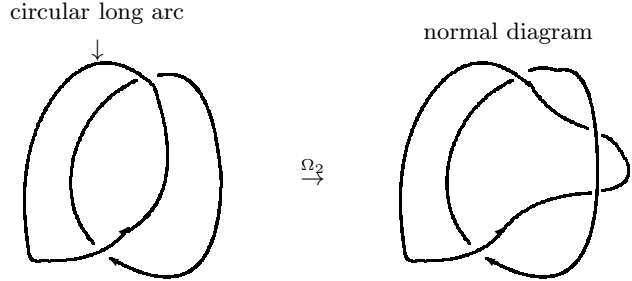

We say that a diagram has {\em proper arcs} if there are no
circles (circular arcs) in the diagram. See Figure 13. We can
choose a diagram $L'$ in such a way that it is divided into long
arcs in a proper way. (A long arc is an arc of the diagram that
may contain virtual crossings.) Such diagrams are called {\em
proper}. It is clear that a proper diagram with $m$ crossings has
$m$ long arcs. We do need that each long arc has two different
final crossing points. For some diagrams this is not true.
However, this can be easily be accomplished by slight deformations
of the diagram, see. Fig. \ref{deform}

\begin{dfn}
An {\em arc} of $L'$ is an oriented interval in the diagram not containing
undercrossings or virtual crossings.
\end{dfn}

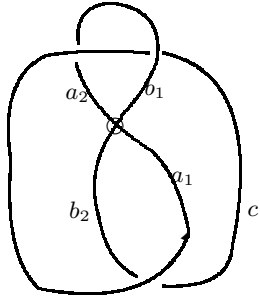
\begin{figure}
\begin{center}
\unitlength 1mm
\linethickness{0.8pt}
\begin{picture}(40.33,58.67)
\put(30.00,17.33){\vector(1,2){0.2}}
\bezier{108}(10.00,10.00)(26.00,7.00)(30.00,17.33)
\bezier{52}(30.00,17.33)(28.67,25.00)(25.00,28.33)
\bezier{64}(25.00,28.33)(16.67,33.67)(15.00,40.00)
\bezier{36}(15.33,42.67)(14.67,47.33)(19.33,48.00)
\bezier{60}(19.33,48.00)(27.00,47.33)(25.67,40.00)
\bezier{36}(25.67,40.00)(24.00,36.00)(21.00,33.00)
\bezier{60}(21.00,32.67)(16.00,26.33)(18.00,19.67)
\bezier{44}(18.00,19.67)(19.00,13.67)(23.00,11.67)
\bezier{52}(26.67,10.33)(35.00,10.00)(35.67,15.00)
\bezier{152}(35.67,15.00)(40.33,38.33)(26.67,41.33)
\bezier{56}(24.67,41.33)(17.00,42.00)(11.00,41.00)
\bezier{72}(11.00,41.00)(4.67,38.00)(6.33,27.00)
\bezier{76}(6.33,27.00)(5.33,14.67)(10.00,10.00)
\put(20.33,31.67){\circle{2.00}}
\put(29.33,24.67){\makebox(0,0)[cc]{$a_1$}}
\put(25.67,36.67){\makebox(0,0)[cc]{$b_1$}}
\put(0.00,158.67){\makebox(0,0)[cc]{$a_2$}}
\put(15.33,35.67){\makebox(0,0)[cc]{$a_2$}}
\put(15.67,20.33){\makebox(0,0)[cc]{$b_2$}}
\put(38.67,20.33){\makebox(0,0)[cc]{$c$}}
\end{picture}
\end{center}
\caption{A knot diagram and its arcs}
\label{knotarcs}
\end{figure}

\begin{ex}
The knot shown in Fig. \ref{knotarcs} has 3 classical crossings,
3 arcs ($a_1$ and $a_2$; $b_1$ and $b_2$; $c$), and 5 virtual arcs
($a_{1},a_{2},b_{1},b_{2},c$).
\end{ex}

The invariant $Q(L)$ is now constructed as follows. Consider all
arcs $a_{i},i=1,\dots, n$ of the diagram $L'$. Consider the set of
formal words $X(L')$ obtained inductively from $a_{i}$ by using
$\circ, \slush, f, f^{-1}$. In order to construct $Q(L')$ we will
factorize $X(L')$ by the equivalence relation generated by the general axioms shown below, plus the
consequences of the specific relations incurred from the diagram.
\bigbreak

Axiomatically, for each $a,b,c\in X(L')$ we identify:

$f^{-1}(f(a))\sim f(f^{-1}(a))\sim a$;

$(a\circ b)\slush b\sim a$;

$(a\slush b)\circ b\sim a$;

$a\sim a\circ a$;

$(a\circ b)\circ c\sim (a\circ c)\circ (b\circ c)$;

$f(a\circ b)\sim f(a)\circ f(b)$.
\bigbreak

With the equivalence generated by these axioms, we get a ``free"
virtual quandle with generators $a_{1},\dots, a_{n}$. We then add
the extra equivalences from the structure of $L':$ \bigbreak

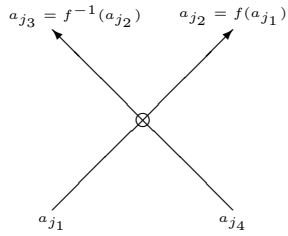
\begin{figure}

\begin{center}
\unitlength 0.6mm
\begin{picture}(50,50)
\put(5,5){\vector(1,1){40}} \put(45,5){\vector(-1,1){40}}
\put(25,25){\circle{3}}
\put(5,2){\makebox(0,0)[cc]{\tiny{$a_{j_1}$}}}
\put(45,2){\makebox(0,0)[cc]{\tiny{$a_{j_4}$}}}
\put(10,48){\makebox(0,0)[cc]{\tiny{$a_{j_3}=f^{-1}({a_{j_2}})$}}}
\put(45,48){\makebox(0,0)[cc]{\tiny{$a_{j_2}=f(a_{j_1})$}}}
\end{picture}
\end{center}

\caption{Relation for a virtual crossing}

\label{virtcross}
\end{figure}

For each classical crossing we write the relation (\ref{clasrel})
just as in the classical case. For each virtual crossing $V$ we
also write relations. Let $a_{j_1},a_{j_2},a_{j_3},a_{j_4}$ be the
four arcs incident to $V$ as it is shown in Fig. \ref{virtcross}.

Then, let us write the relations:

\begin{equation} a_{j_2}= f(a_{j_1}) \label{former2}\end{equation}

and

\begin{equation} a_{j_3}=f(a_{j_4}) \label{former3}\end{equation}

So, the virtual quandle $Q(L)$ is the quandle, generated by all
arcs $a_{i},i=1,\dots, n$, all linear relations at classical
vertices and all relations (\ref{former2},\ref{former3}) at
virtual vertices.

\begin{th}\label{INV}
The quandle $Q(L)$ is a link invariant.
\end{th}

\begin{proof}

First, note that two proper diagrams generate
isotopic virtual links if and only if one can be deformed to the other
by using a sequence of virtual Reidemeister. Indeed,
if while isotopy a circular long link occurs, we
can modify the isotopy by applying the first classical
Reidemeister move to this long arc and subdividing
it into two parts.

We have to show that $Q(L)$ is invariant under virtual
Reidemeister moves.

The invariance of $Q(L)$ under all classical moves is checked
in the same way as that of the ordinary quandle.

Let us check now the invariance of $Q$ under purely virtual
Reidemeister moves.

The first virtual Reidemeister move is shown in Fig. \ref{1move}.
In the initial local picture we have one local
generator $a$.
Here we just add a new generator $b$ and two
coinciding relations: $b=f(a)$. Thus, it does
not change the virtual quandle at all.

\begin{figure}
\unitlength 0.4mm
\begin{center}
\begin{picture}(100,55)
\thicklines
\put(10,10){\line(0,1){30}}
\put(40,40){\vector(0,-1){30}}
\put(25,40){\oval(30,20)[t]}
\put(25,45){\makebox(0,0)[cc]{$a$}}

\put(50,25){\makebox(0,0)[cc]{$\longrightarrow$}}

\put(60,10){\line(1,1){30}}
\put(60,40){\vector(1,-1){30}}
\put(75,40){\oval(30,20)[t]}
\put(75,25){\circle{4}}
\put(55,5){\makebox(0,0)[cc]{$a$}}
\put(75,45){\makebox(0,0)[cc]{$b$}}
\put(95,5){\makebox(0,0)[cc]{$a$}}

\end{picture}

\end{center}
\caption{Invariance of $Q$ under the $1$--st virtual move}
\label{1move}
\end{figure}
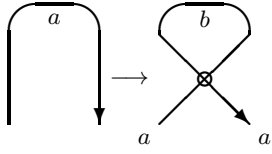

The case of inverse orientation at the crossings
gives us $b=f^{-1}(a)$ that does not change the situation.

For each next relation, we will check only one case
of arc orientation.

The second move (see Fig. \ref{2move}) adds two generators
$c$ and $d$ and two pairs of coinciding relations:
$c=f(a),d=f^{-1}(b)$. Thus, the quandle $Q$ stays the
same.

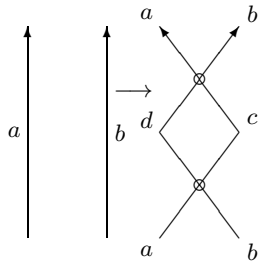
\begin{figure}
\begin{center}
\begin{picture}(100,100)

\put(10,10){\vector(0,1){80}}
\put(40,10){\vector(0,1){80}}
\put(5,50){\makebox(0,0)[cc]{$a$}}
\put(45,50){\makebox(0,0)[cc]{$b$}}

\put(50,65){\makebox(0,0)[cc]{$\longrightarrow$}}

\put(60,10){\line(3,4){30}}
\put(90,10){\line(-3,4){30}}
\put(60,50){\vector(3,4){30}}
\put(90,50){\vector(-3,4){30}}
\put(75,30){\circle{4}}
\put(75,70){\circle{4}}
\put(55,5){\makebox(0,0)[cc]{$a$}}
\put(95,5){\makebox(0,0)[cc]{$b$}}
\put(55,55){\makebox(0,0)[cc]{$d$}}
\put(95,55){\makebox(0,0)[cc]{$c$}}
\put(55,95){\makebox(0,0)[cc]{$a$}}
\put(95,95){\makebox(0,0)[cc]{$b$}}

\end{picture}
\end{center}

\caption{Invariance of $Q$ under the $2$--nd virtual move}
\label{2move}

\end{figure}

In the case of the third Reidemeister move we have six ``exterior
arcs'': three incoming $(a,b,c)$ and three outgoing $(p,q,r)$, see
Fig. \ref{3move}. In both cases we have
$p=f^{2}(a),q=b,r=f^{-2}(c)$. The three interior arcs are
expressed in $a,b,c$, and give no other relations.

\begin{figure}
\begin{center}
\begin{picture}(200,100)
\put(50,10){\vector(0,1){80}}
\put(10,50){\vector(1,0){80}}
\put(10,10){\line(2,1){40}}
\put(50,30){\line(1,1){20}}
\put(70,50){\vector(1,2){20}}

\put(50,50){\circle{4}}
\put(50,30){\circle{4}}
\put(70,50){\circle{4}}

\put(5,50){\makebox(0,0)[cc]{$a$}}
\put(5,5){\makebox(0,0)[cc]{$b$}}
\put(50,5){\makebox(0,0)[cc]{$c$}}

\put(95,50){\makebox(0,0)[cc]{$p$}}
\put(95,95){\makebox(0,0)[cc]{$q$}}
\put(50,95){\makebox(0,0)[cc]{$r$}}

\put(100,60){\makebox(0,0)[cc]{$\longrightarrow$}}

\put(150,10){\vector(0,1){80}}
\put(110,50){\vector(1,0){80}}
\put(110,10){\line(1,2){20}}
\put(130,50){\line(1,1){20}}
\put(150,70){\vector(2,1){40}}
\put(150,50){\circle{4}}
\put(130,50){\circle{4}}
\put(150,70){\circle{4}}

\put(105,50){\makebox(0,0)[cc]{$a$}}
\put(105,5){\makebox(0,0)[cc]{$b$}}
\put(150,5){\makebox(0,0)[cc]{$c$}}

\put(195,50){\makebox(0,0)[cc]{$p$}}
\put(195,95){\makebox(0,0)[cc]{$q$}}
\put(150,95){\makebox(0,0)[cc]{$r$}}

\end{picture}
\end{center}

\caption{Invariance of $Q$ under the $3$--rd virtual move}
\label{3move}
\end{figure}
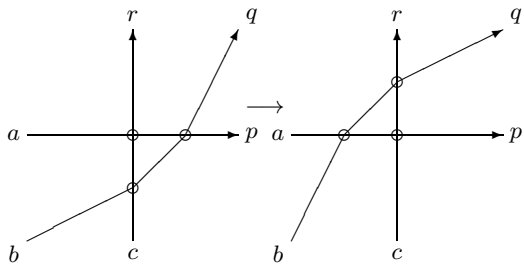

Finally, let us check the mixed move. We will check
the only version of it, see Fig. \ref{mixed}

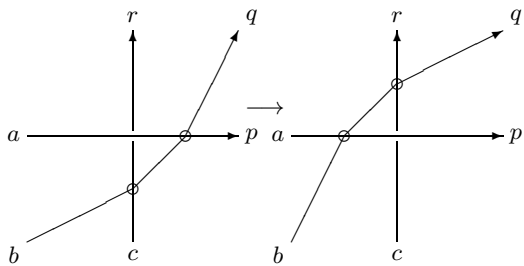
\begin{figure}
\begin{center}

\begin{picture}(200,100)
\put(50,10){\line(0,1){38}}
\put(50,52){\vector(0,1){38}}
\put(10,50){\vector(1,0){80}}
\put(10,10){\line(2,1){40}}
\put(50,30){\line(1,1){20}}
\put(70,50){\vector(1,2){20}}

\put(50,30){\circle{4}}
\put(70,50){\circle{4}}

\put(5,50){\makebox(0,0)[cc]{$a$}}
\put(5,5){\makebox(0,0)[cc]{$b$}}
\put(50,5){\makebox(0,0)[cc]{$c$}}

\put(95,50){\makebox(0,0)[cc]{$p$}}
\put(95,95){\makebox(0,0)[cc]{$q$}}
\put(50,95){\makebox(0,0)[cc]{$r$}}

\put(100,60){\makebox(0,0)[cc]{$\longrightarrow$}}

\put(150,10){\line(0,1){38}}
\put(150,52){\vector(0,1){38}}
\put(110,50){\vector(1,0){80}}
\put(110,10){\line(1,2){20}}
\put(130,50){\line(1,1){20}}
\put(150,70){\vector(2,1){40}}

\put(130,50){\circle{4}}
\put(150,70){\circle{4}}

\put(105,50){\makebox(0,0)[cc]{$a$}}
\put(105,5){\makebox(0,0)[cc]{$b$}}
\put(150,5){\makebox(0,0)[cc]{$c$}}

\put(195,50){\makebox(0,0)[cc]{$p$}}
\put(195,95){\makebox(0,0)[cc]{$q$}}
\put(150,95){\makebox(0,0)[cc]{$r$}}

\end{picture}

\end{center}
\caption{Invariance of $Q$ under the mixed move}
\label{mixed}
\end{figure}

In both pictures
we have three incoming edges $a,b,c$ and three outgoing
edges $p,q,r$. In the first case we have relations:
$p=f(a),q=b,r=f^{-1}(c)\circ a$. In the
second case we have: $p=f(a), q=b,
r=f^{-1}(c\circ f(a))$.

The distributivity relation $f(x\circ y)=f(x)\circ f(y)$
implies $f^{-1}(c)\circ a=f^{-1}(c\circ f(a))$.
Hence, two virtual quandles
before the mixed move and after the mixed move coincide.

The other cases of the mixed move lead to other
relations all equivalent to $f(x\circ y)=f(x)\circ f(y)$.

This completes the proof of the theorem.

\end{proof}

This leads to a construction of a polynomial invariant (the
``easier'' part of the paper \cite{Ma3}). First, we take $a\circ
b$ to be $ta+(1-t)b$ and $f(a)$ to be $sa$, where $t$ and $s$ are
independent commuting variables. After this, the normalized
determinant of this matrix gives an invariant.

It was proved recently by R.Fenn (considering virtual knots as
closures of virtual braids, this idea was inspirated by H.Morton)
that this polynomial in fact coincides with that proposed by
Sawollek and Silver-Williams and that coming from the Alexander
biquandle (up to a variable change) as discussed by Kauffman and
Radford. Fenn's method also proves that we get nothing new if we
try to use the linear biquandle structure at classical crossing
and linear automorphism $f$ at virtual crossings.

In the present paper, we will construct a more general invariant
in section \ref{sec} (which works for the case of colored links).

\subsection{Jones polynomial for virtual knots \\ and involutory
quandles}

We use a generalization of the bracket state summation model for
the Jones polynomial to extend it to virtual knots and links.  We
call a diagram in the plane {\em purely virtual} if the only
crossings in the diagram are virtual crossings. Each purely
virtual diagram is equaivalent by the virtual moves to a disjoint
collection of circles in the plane.

\begin{figure}
\centering\includegraphics[width=300pt]{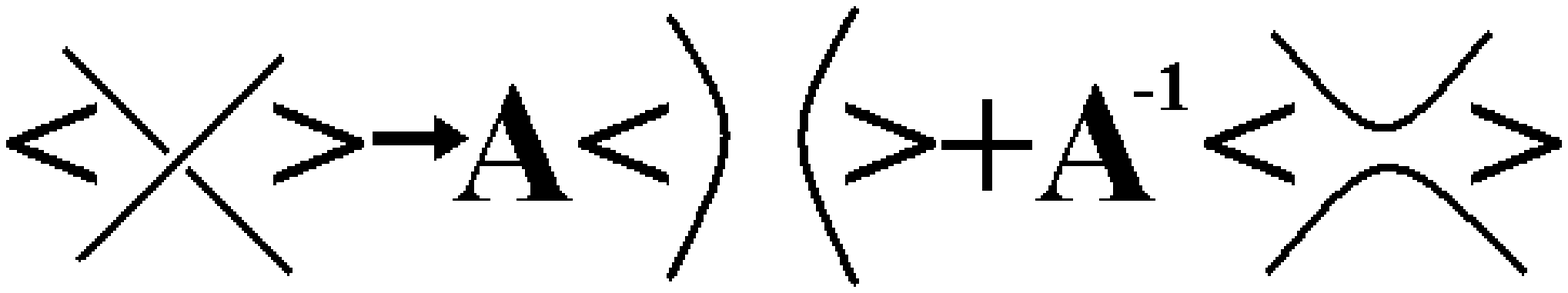} \label{kf}
\caption{}
\end{figure}

Given a link diagram $K$, a state $S$ of this diagram is obtained
by choosing a smoothing for each crossing in the diagram and
labelling that smoothing with either $A$ or $A^{-1}$ according to
the convention that a counterclockwise rotation of the
overcrossing line sweeps two regions labelled $A$, and that a
smoothing that connects the $A$ regions is labelled by the letter
$A$, or, diagrammatically, see Fig. \ref{kf}

 Then, given a state $S$, one has the
evaluation $\langle K|S\rangle$ equal to the product of the labels
at the smoothings, and one has the evaluation $||S||$ equal to the
number of loops in the state (the smoothings produce purely
virtual diagrams).  One then has the formula

\begin{equation}
\langle K\rangle  = \Sigma_{S}\langle K|S\rangle
d^{||S||-1},\label{frm1}
\end{equation} where the summation runs over the states $S$ of the
diagram $K$, and $d = -A^{2} - A^{-2}.$ This state summation is
invariant under all classical and virtual moves except the first
Reidemeister move. The bracket polynomial is normalized to an
invariant $f_{K}(A)$ of all the moves by the formula  $f_{K}(A) =
(-A^{3})^{-w(K)}\langle K\rangle$ where $w(K)$ is the writhe of
the (now) oriented diagram $K$. The writhe is the sum of the
orientation signs ($\pm 1)$ of the crossings of the diagram. The
Jones polynomial, $V_{K}(t)$ is given in terms of this model by
the formula \begin{equation} V_{K}(t) =
f_{K}(t^{-1/4}).\label{frm2}\end{equation}

The reader should note that this definition is a direct
generalization to the virtual category of the state sum model for
the original Jones polynomial. It is straightforward to verify the
invariances stated above. In this way one has the Jones polynomial
for virtual knots and links, see \cite{KaV}.

In terms of the interpretation of virtual knots as stabilized
classes of embeddings of circles into thickened surfaces, our
definition coincides with the simplest version of the Jones
polynomial for links in thickened surfaces. In that version one
counts all the loops in a state the same way, with no regard for
their isotopy class in the surface. It is this equal treatment
that makes the invariance under handle stabilization work. With
this generalized version of the Jones polynomial, one has again
the problem of finding a geometric/topological interpretation of
this invariant. There is no fully satisfactory topological
interpretation of the original Jones polynomial and the problem is
inherited by this generalization.  \bigbreak

In \cite{KaV2}, the following theorem was proved.

\begin{th} To
each non-trivial classical knot diagram of one component $K$ there
is a corresponding  non-trivial virtual knot diagram $Virt(K)$
with unit Jones polynomial.
\end{th}

This Theorem is a key ingredient in the problems involving virtual
knots. Here is a sketch of its proof. The proof uses two
invariants of classical knots and links that generalize to
arbitrary virtual knots and links. These invariants are the {\em
Jones polynomial} and the {\em involutory quandle} denoted by the
notation $IQ(K)$ for a knot or link $K.$

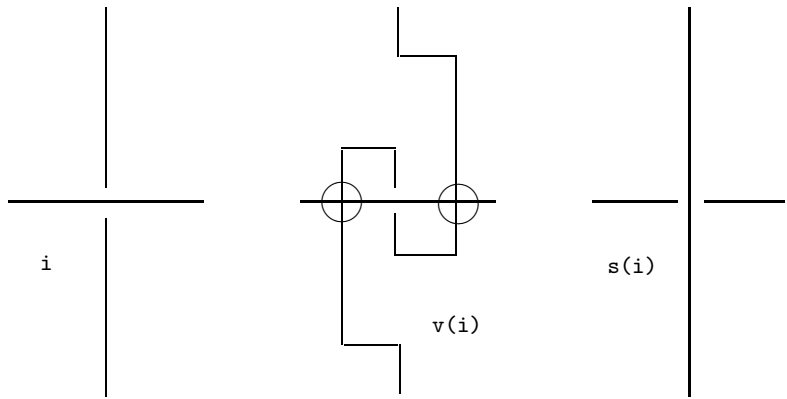
\begin{figure}
\begin{center}
{\tt    \setlength{\unitlength}{0.92pt}
\begin{picture}(326,162)
\thinlines    \put(190,141){\line(0,-1){82}}
              \put(242,34){\makebox(41,41){s(i)}}
              \put(170,10){\makebox(41,41){v(i)}}
              \put(1,35){\makebox(41,42){i}}
              \put(292,81){\line(1,0){33}}
              \put(246,81){\line(1,0){35}}
              \put(286,161){\line(0,-1){160}}
              \put(143,81){\circle{16}}
              \put(191,80){\circle{16}}
              \put(190,59){\line(-1,0){25}}
              \put(167,22){\line(0,-1){20}}
              \put(144,22){\line(1,0){22}}
              \put(143,102){\line(0,-1){80}}
              \put(165,103){\line(-1,0){22}}
              \put(165,87){\line(0,1){15}}
              \put(165,59){\line(0,1){17}}
              \put(167,141){\line(1,0){23}}
              \put(166,161){\line(0,-1){20}}
              \put(126,81){\line(1,0){80}}
              \put(46,74){\line(0,-1){73}}
              \put(46,161){\line(0,-1){74}}
              \put(6,81){\line(1,0){80}}
\end{picture}}
\end{center}

 \caption{Switching and Virtualizing a Crossing}
\label{fj3}
\end{figure}

Given a crossing $i$ in a link diagram, we define $s(i)$ to be the
result of {\em switching} that crossing so that the undercrossing
arc becomes an overcrossing arc and vice versa. We also define the
{\em virtualization} $v(i)$ of the crossing by the local
replacement indicated in Figure \ref{fj3}. In this Figure we
illustrate how in the virtualization of the crossing the  original
crossing is replaced by a crossing that is flanked by two virtual
crossings. \bigbreak

Suppose that $K$ is a (virtual or classical) diagram with a
classical crossing labelled $i.$  Let $K^{v(i)}$ be the diagram
obtained from $K$ by virtualizing the crossing $i$ while leaving
the rest of the diagram just as before. Let $K^{s(i)}$ be the
diagram obtained from $K$ by switching the crossing $i$ while
leaving the rest of the diagram just as before. Then it follows
directly from the definition of the Jones polynomial that
$$V_{K^{s(i)}}(t) = V_{K^{v(i)}}(t).$$ \noindent As far as the
Jones polynomial is concerned, switching a crossing and
virtualizing a crossing look the same. \bigbreak

\begin{figure}
\begin{center}
{\tt    \setlength{\unitlength}{0.92pt}
\begin{picture}(214,164)
\thinlines    \put(197,64){\makebox(16,16){b}}
              \put(136,144){\makebox(22,18){c =}}
              \put(17,143){\makebox(22,18){c =}}
              \put(165,143){\makebox(19,20){ab}}
              \put(113,64){\makebox(16,16){b}}
              \put(67,64){\makebox(16,16){b}}
              \put(43,1){\makebox(17,19){a}}
              \put(47,142){\makebox(19,20){ab}}
              \put(1,63){\makebox(16,16){b}}
              \put(165,1){\makebox(17,19){a}}
              \put(185,142){\line(0,-1){82}}
              \put(138,82){\circle{16}}
              \put(186,81){\circle{16}}
              \put(185,60){\line(-1,0){25}}
              \put(162,23){\line(0,-1){20}}
              \put(139,23){\line(1,0){22}}
              \put(138,103){\line(0,-1){80}}
              \put(160,104){\line(-1,0){22}}
              \put(160,88){\line(0,1){15}}
              \put(160,60){\line(0,1){17}}
              \put(162,142){\line(1,0){23}}
              \put(161,162){\line(0,-1){20}}
              \put(121,82){\line(1,0){80}}
              \put(41,75){\line(0,-1){73}}
              \put(41,162){\line(0,-1){74}}
              \put(1,82){\line(1,0){80}}
\end{picture}}
\end{center}
\caption{$IQ(Virt(K)) = IQ(K)$} \label{fj4}
\end{figure}
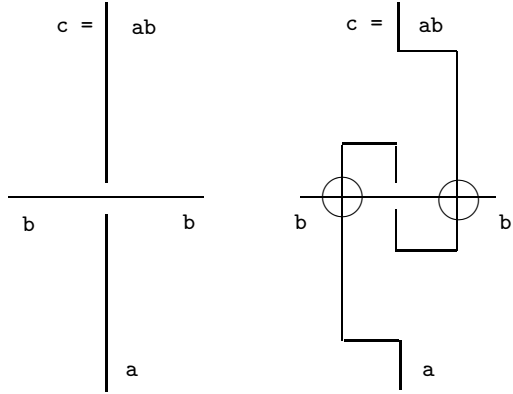

The involutory quandle \cite{KNOTS} is an algebraic invariant
equivalent to the fundamental group of the double branched cover
of a knot or link in the classical case. In this algebraic system
one associates a generator of the algebra $IQ(K)$ to each arc of
the diagram $K$ and there is a relation of the form $c = ab$ at
each crossing, where $ab$ denotes the (non-associative) algebra
product of $a$ and $b$ in $IQ(K).$ See Figure \ref{fj4}. In this
Figure we have illustrated through the local relations the fact
that
$$IQ(K^{v(i)}) = IQ(K).$$ \noindent As far the involutory quandle
is concerned, the original crossing and the virtualized crossing
look the same. \bigbreak

If a classical knot is actually knotted, then its involutory
quandle is non-trivial \cite{W}. Hence if we start with a
non-trivial classical knot, we can virtualize any subset of its
crossings to obtain a virtual knot that is still non-trivial.
There is a subset $A$ of the crossings of a classical knot $K$
such that the knot $SK$ obtained by switching these crossings is
an unknot.  Let $Virt(K)$ denote the virtual diagram obtained from
$A$ by virtualizing the crossings in the subset $A.$  By the above
discussion the Jones polynomial of $Virt(K)$ is the same as the
Jones polynomial of $SK$, and this is $1$ since $SK$ is unknotted.
On the other hand, the $IQ$ of $Virt(K)$ is the same as the $IQ$
of $K$, and hence if $K$ is knotted, then so is $Virt(K).$   We
have shown that $Virt(K)$ is a non-trivial virtual knot with unit
Jones polynomial.  This completes the proof of the Theorem.

If there exists a classical knot with unit Jones polynomial, then
one of the knots $Virt(K)$ produced by this Theorem may be
equivalent to  a classical knot.  It is an intricate task to
verify that specific examples of $Virt(K)$ are not classical; very
special partial case was considered in \cite{SW2}. This has led to
an investigation of new invariants for virtual knots. In this
investigation a number of issues appear. One can examine the
combinatorial generalization of the fundamental group (or quandle)
of the virtual knot and sometimes one can prove by pure algebra
that the resulting group is not classical. This is related to
observations by Silver and Williams \cite{SW}, Manturov
\cite{Ma6,Ma8} and by Satoh \cite{Satoh} showing that the
fundamental group of a virtual knot can be interpreted as the
fundamental group of the complement of a torus embedded in four
dimensional Euclidean space. A very fruitful line of new
invariants comes about by examining the biquandles and virtual
quandles. Flat virtual diagrams are seldom trivial. If we can
verify that the flat knot $F(Virt(K))$ is non-trivial, then
$Virt(K)$ is non-classical. In this way the search for classical
knots with unit Jones polynomial expands to the exploration of the
structure of the infinite collection of virtual knots with unit
Jones polynomial, for a detailed description of this problem, see
\cite{FKM}.

Another way of putting this Theorem is as follows: In the arena of
knots in thickened surfaces there are many examples of knots with
unit Jones polynomial. Might one of these be equivalent via handle
stabilization to a classical knot? In \cite{KUP} Kuperberg shows
the uniqueness of the  embedding of minimal genus in the stable
class for a given virtual link. The minimal embedding genus can be
strictly less than the number of virtual crossings in a diagram
for the link.  There are many problems associated with this
phenomenon.

There are generalizations of the Jones polynomial that involve the
surface representation of virtual knots. To begin with, one can keep track of the isotopy classes of the
curves in the state expansion of the bracket polynomial for a knot embedded in a surface. This gives a
surface bracket polynomial \cite{KaV4} that can be used in tandem with Kuperberg's results \cite{KUP}
to determine the minimal surface embedding genus for some virtual knots and links. In \cite{KaV4} this
method is used to show that the Kishino diagram has genus two. Another approach \cite{Ma8}
uses a relative of the surface bracket polynomial, and includes in the equivalence relation of the curves in
the states, the stabilization of the surfaces themselves. In this way the Manturov invariant is
an  becomes an element of a module over the Laurent polynomial ring in one variable, and
is strictly stronger than the original extension of the Jones polynomial to
virtuals. We do not yet know the relative strengths of these two methods.
\bigbreak

\subsection{A Common Construction --- \\ The Virtual Biquandle}

Consider an oriented virtual link diagram $L$ and divide it into
arcs. Write down the biquandle relations at each classical
crossing, and write down the virtual quandle relation at each
virtual crossing.

The full list of operations and axioms for this object, {\em the virtual biquandle}
(making this related algebra an invariant of virtual knots and links) is as follows:

\begin{enumerate}

\item We have operations $\UR{}\;,\LR{}\;,\UL{}\;,\LL{}$ as above
and the invertible operation $f$ acting on the whole algebra.

\item Biquandle operations satisfy all conditions described previously.

\item The operation $f$ ia a biquandle automorphism. Thus, for
instance, $f(a\UR{b})=f(a)\UR{f(b)}$; the same statements are true
for for $\LR{~},\UL{~},\LL{~}$.

\end{enumerate}

The object defined in this manner is called {\em the virtual
biquandle of the link}. Now, we consider a larger context,
admitting more general constructions at virtual crossings.

Let us now give the definition of the {\em formal virtual
biquandle} as some algebraic object satisfying certain axioms; in
the same lines, we shall define the virtual biquandle of a link by
interpreting the biquandle operations according to the crossing
structure.

\begin{figure}
\unitlength 1.8pt
\begin{center}
\begin{picture}(50,50)
\put(5,5){\vector(1,1){40}} \put(45,5){\vector(-1,1){40}}
\put(25,25){\circle{3}} \put(0,0){a} \put(50,0){b}
\put(50,50){a.b}\put(0,50){b|a}
\end{picture}
\end{center}
\caption{Binary relation at a virtual crossing} \label{bll}
\end{figure}
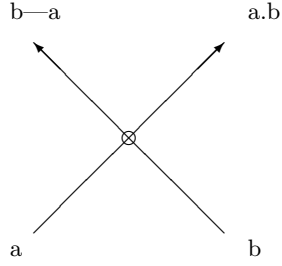

Let us now write down the axioms in the general case, when virtual
crossings are endowed with binary operations, $a.b$ and $b|a$, see
Fig. \ref{bll}.

The axioms for the classical case are the same as above.

The axioms for purely virtual moves are:

\begin{enumerate}

\item The first Reidemeister move, see Fig. \ref{FRM}

\begin{equation}\forall a\; \exists x: \left\{\begin{array}{c}x.a=a \cr
a|x=x\end{array}\right.
\end{equation}

\begin{equation} \forall a\; \exists y: \left\{\begin{array}{c}a.y=y \cr
y|a=a\end{array}\right.
\end{equation}

\begin{figure}
\centering\includegraphics[width=200pt]{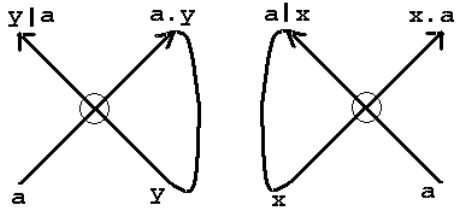}
\caption{Relations coming from the first virtual move} \label{FRM}
\end{figure}

\item The second oriented Reidemeister move, see Fig.
\ref{SORM}.1.

\begin{equation} \forall a,b: (a.b)|(b|a)=a\mbox{ and }
(b|a).(a.b)=b\end{equation}

The second unoriented Reidemeister moves, see Fig. \ref{SORM}.2,
\ref{SORM}.3.

\begin{equation}
\left\{\begin{array}{c} \forall a,b:\exists! x:
b|(a|x)=x;(a|x).b=a;\cr \forall a,b:\exists! y: (y|b)=a;
(b.y)|a=b.
\end{array} \right.
\end{equation}

\begin{figure}
\centering\includegraphics[width=300pt]{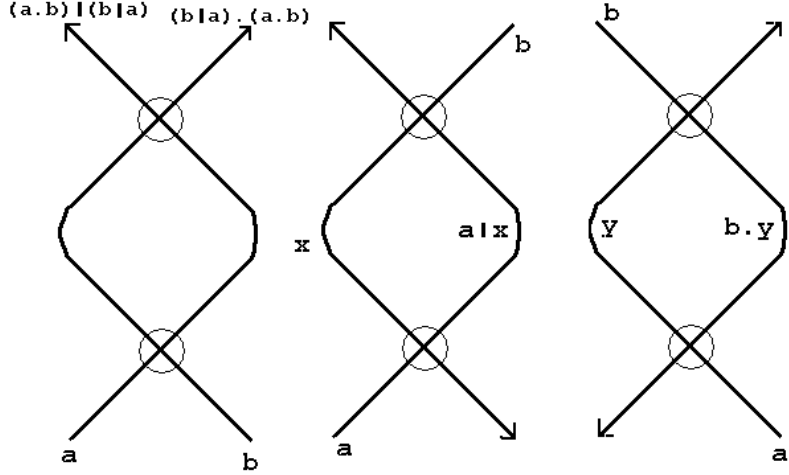}
\caption{Relations coming from the second virtual move}
\label{SORM}
\end{figure}

\item The third Reidemeister move, see Fig. \ref{TRM} (upper and
lower pictures).

The two cases are:

\begin{equation}
\forall a,b,c: \left\{\begin{array}{c} (a.b).c=(a.(c|b)).(b.c)\cr
(b|a).(c|(a.b))=(b.c)|(a.(c|b)) \cr (c|(a.b))|(b|a)=(c|b)|a
\end{array} \right.
\end{equation}

and

\begin{equation}
\forall a,b,c: \exists! x,y: \left\{\begin{array}{c}
(c|(b|(a|x)))=x;\cr b.(c.(a.y))=y; \cr (a|x).b=(a.y)|c;\cr
y|a=(b|(a|x)).c;\cr x.a=(c.(a.y))|b.\end{array}\right.
\end{equation}

\begin{figure}
\centering\includegraphics[width=300pt]{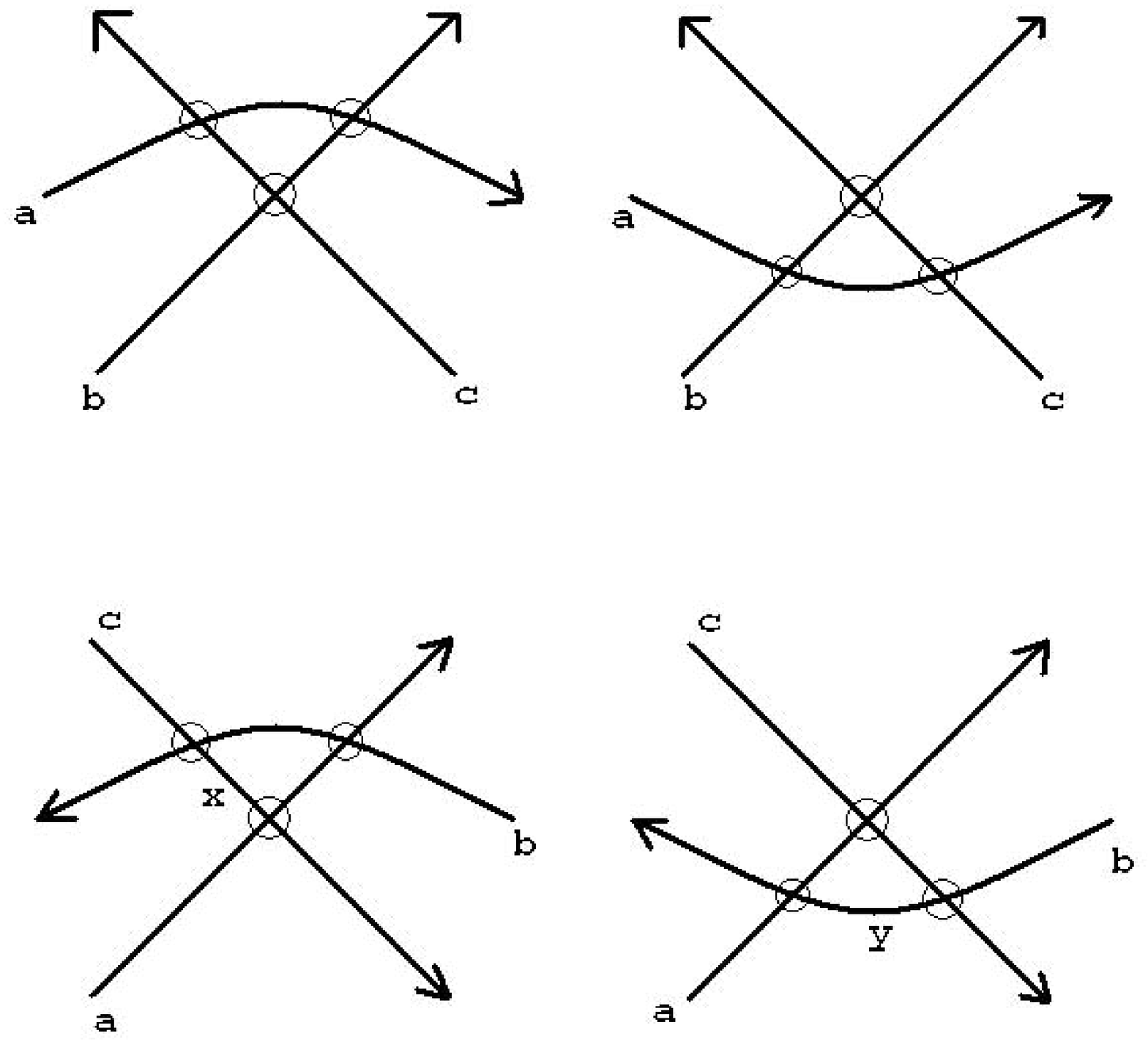}
\caption{Relations coming from the third virtual move} \label{TRM}
\end{figure}

\end{enumerate}

\begin{re}
In Figs. \ref{SORM} and \ref{TRM} we label only some arcs; all
other arcs can be labelled according to the rule shown in Fig.
\ref{bll}.
\end{re}

\begin{figure}
\centering\includegraphics[width=300pt]{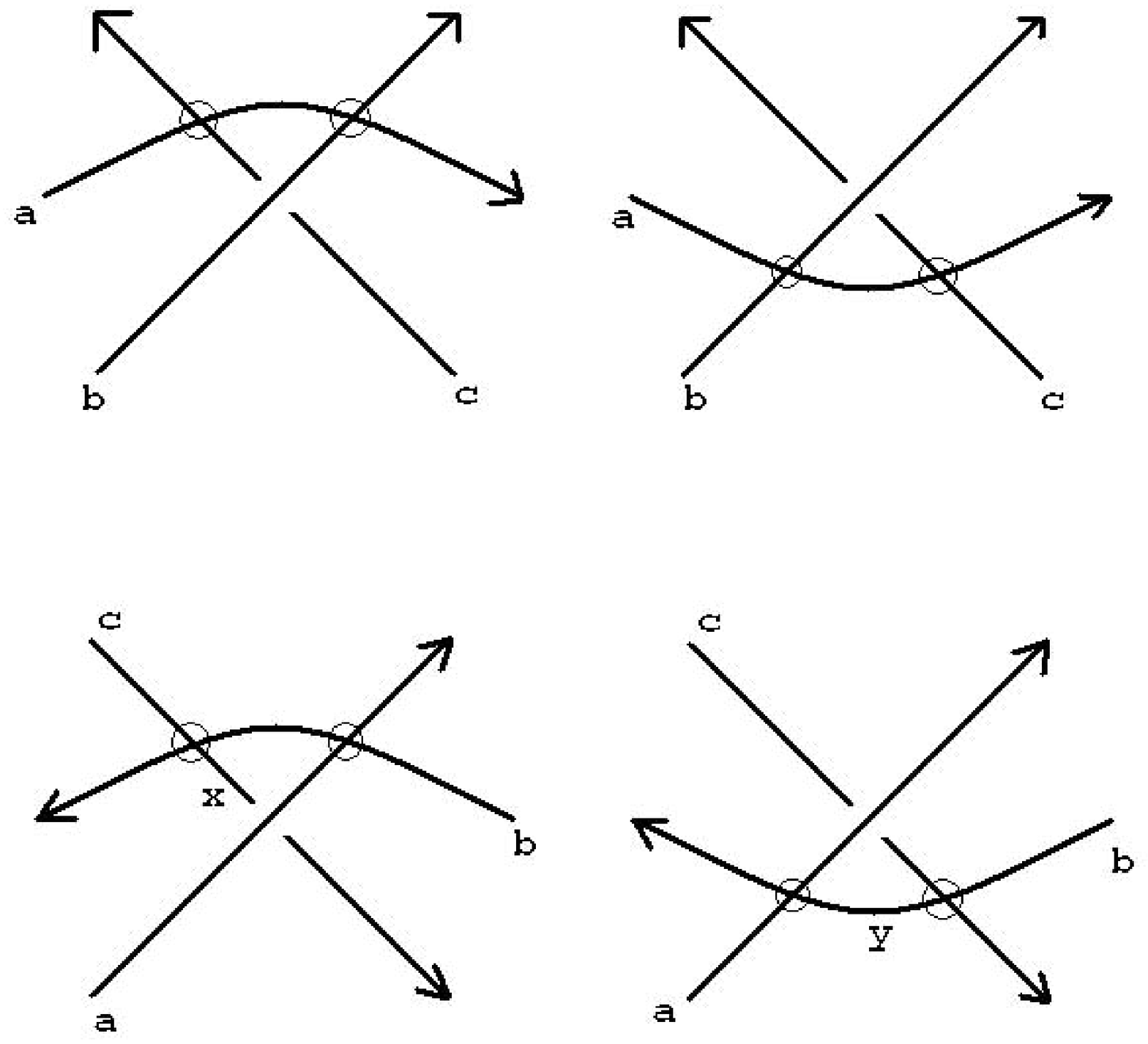} \caption{Relations
coming from the mixed move} \label{MM}
\end{figure}

For the semivirtual move we have the following pictures in Fig.
\ref{MM} and relations:

\begin{equation}
\forall a,b,c: \left\{\begin{array}{c}
(a\UR{b})|c=(a|(c.b))\UR{b|c};\cr
(b\LR{a})|(c.(a\UR{b})=(b|c)\LR{a|(c.b)};\cr
c.(a\UR{b}).(b\LR{a})=(c.b).a
\end{array}\right.
\end{equation}

and

\begin{equation}
\forall a,b,c: \exists! x,y:\left\{\begin{array}{c}
(b|(a\LL{x})).c=x ;\cr a.(b.(c\LR{y}))=y;\cr
(a\LL{x}).b=y\LL{c};\cr c|(b|(a\LL{x}))=(b.(c\UL{y}))|a ;\cr
x\UL{a}=(c\UL{y})|b
\end{array}\right.
\end{equation}

The proof is trivial by design. It is left to the reader.

This common structure should have a lot of realizations. We
indicate one of them in the following section.

\section{Presentations}

In this section, we will describe several explicit
constructions of the virtual quandle and biquandle invariants described above, using specific
algebraic structures for the representation.

\subsection{A presentation for virtual biquandles}

Virtual biquandles admit the following linear representation: at
classical crossings we have the Alexander biquandle operations
(with generators $t,s$) and set $(a.b)=(1+q b)a; (b'a)=(1-q a)b$
for some new generator $q$ that commutes with $t,s$ such that
$q^{2}=(t-1)q=(s-1)q=0$. We are going to investigate other models
of the virtual biquandle with binary operations at virtual
crossings.

We shall discuss other presentations of the virtual biquandle
construction in our further publications. Later in this section,
we deal with long quandles, virtual quandles, and some
generalizations of them.

\subsection{Formal Power Functions and Conjugation}

For the virtual quandles, one can take a group together with the
following two possible operations:

\begin{enumerate}

\item $a\circ b= b^{n}ab^{-n},$ $f(a)=qaq^{-1}$, where $q$ is a
new generator for this group, and $n$ is an integer.

\item $a\circ b=ba^{-1}b$, $f(a)=qaq^{-1}$ or $f(a)=qa^{-1}q$,
where $q$ is a new generator for the group.

\end{enumerate}

In this way, we associate a group to each virtual knot in two possible ways.

One can also associate a finite group $G$ with the virtual quandle
operations defined as above, and look at the set of all maps
$\Gamma(L)\to G$, where $\Gamma(L)$ is the virtual knot quandle,
and $G$ is a finite virtual quandle. The set of such mappings is
finite for any finite group (we can fix the images of generators
of the quandle thus defining the map completely), so the number of
``colorings'', i.e., mappings, is finite and invariant under
generalised Reidemeister moves, for more details see \cite{Ma2}.

These ideas can also be used for the case of biquandles and virtual
biquandles.

More precisely, one can say the following:  if we have an abstract
knot invariant ${\cal G}$  constructed according to some axioms
(quandle, virtual quandle or biquandle) and a finite object $G$
satisfying the same axioms, then the cardinality of the set of
mappings ${\cal G}(L)\to G$ is a knot invariant.

Thus there is a way to search for finite--valued invariants (colorings)
by finding {\bf finite} sets satisfying such axioms.

\subsection{The Quaternionic biquandle}

It turns out that there are beautiful, explicit biquandle
constructions that give powerful results. One class of such
constructions come from the quaternions. This idea and the
following formulae are due to R.A. Fenn and A. Bartholomew
\cite{BF}.

The quaternionic biquandle is defined by the following operations
where $i^2 = j^2 = k^{2} = ijk = -1, ij = -ji = k,jk = -kj = i, ki
= -ik =j$ in the associative, non-commutative algebra of the
quaternions. The elements $a,b, \cdots$ are in a module over the
ring of integer quaternions.

$$a \UR{b}=j\cdot a+ (1+i)\cdot b,$$

$$a \LR{b}=-j\cdot a+ (1+i)\cdot b,$$

$$a \UL{b}=j\cdot a+ (1-i) \cdot b,$$

$$a \LL{b}=-j\cdot a+ (1-i)\cdot b.$$

Amazingly, one can verify that these operations satisfy the axioms for the biquandle.

For more about quaternionic and other non-commutative buquandles,
see \cite{BuF}. \bigbreak

\begin{figure}
\centering\includegraphics[width=180pt]{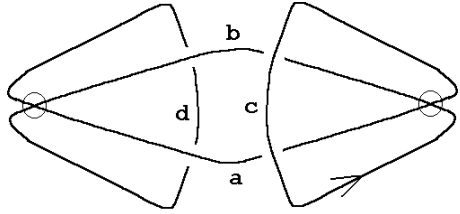} \caption{The
Kishino knot} \label{Ksl}
\end{figure}

Let us look at what happens with the Kishino knot (see labeling
shown in Fig. \ref{Ksl}.

Let us first consider the rightmost crossings (that deal with arcs
a,b,c).

We have:

$$c\LR{a}=-j\cdot c+(1+i)a; a\UR{c}=j\cdot a+(1+i)\cdot c.$$

Now,

$$a\UR{c} \, \LL{c \, \LR{a}} = c \, \mbox{ and } c \, \LR{a} \, \UL{a \, \UR{c}} = b$$

give us

\begin{equation}
b=j(-jc+(1+i)a)+(1-i)(ja+(1+i)c)=3c+2(j-k)a;\label{form4}\end{equation}

\begin{equation}
c=-j(ja+(1+i)c)+(1-i)(-jc+(1+i)a)=3a-2(j-k)c.\label{form5}\end{equation}

It follows that $3a=3b$.

Now, let us write the relations for the left crossings. We have:

$$b\LL{d}=-jb+(1-i)d; d\UL{b}=jd+(1-i)b.$$

Now,

$$d \, \UL{b} \, \LR{b \, \LL{d}} = a \, \mbox{ and } b \, \LL{d} \, \UR{d \, \UL{b}} = d.$$

So,

\begin{equation}a=
-j(jd+(1-i)b)+(1+i)(-jb+(1-i)d)=3d-2(j+k)b\label{form6}\end{equation}

\begin{equation}d= j (-jb+(1-i)d)
+(1+i)(jd+(1-i)b)=3b+2(j+k)d\label{form7}\end{equation}

From these equations we can also conclude that $3a=3b$. Let us
show that the module we get is not just the module generated by
$a$ where $a=b$.

Indeed, by tensoring the above equations with ${\bf Z}_{3}$, we
get: $b=2(j-k)a, c=-2(j-k)c;a=-2(j+k)b; d=2(j+k)d$. Thus, $b$ is
out of picture, and we get:

$c=-2(j-k)c; a=2 i a; d=2(j+k)d$.

But this implies $(1-2i)a=0$ from which we deduce
$2a=5a=(1+2i)(1-2i)a=0$, whence $a=0$.

Thus we are left with $c=2(j-k)c, d=2(j+k)d$ and this is certainly
a non-trivial module over ${\bf Z}_{3}$.

 Thus, this linear biquandle is not isomorphic to the
free one--dimensional linear space over integer quaternions. From
this we conclude that the Kishino diagram is a non-trivial virtual
knot. The proof given above is a simplification of the method used
by Fenn and Bartholomew. In any case, this proof is perhaps the
most direct verification for the detection of the Kishino diagram!
\bigbreak

\section{Infinite-dimensional Lie algebras}

 Here we propose one new method to obtain quandles,
biquandles, virtual quandles, etc, see \cite{Ma4}. For simplicity,
we will deal with quandles. We know that quandles can be well
defined on discrete groups by putting, say, $a\circ b=bab^{-1}$.
The question is: is it possible to arrange a quandle operation
which would act on a geometric group (Lie group)? Well, each of
the relations

\begin{equation}
\left\{\begin{array}{ccc} a\circ b &=& bab^{-1}\cr a\circ b
&=&b^{n}a b^{-n} \cr a\circ b & = & ba^{-1} b
\end{array}\right.\label{qrela}
\end{equation}
works well in \underline{any} group (i.e. satisfies all the
quandle axioms), but what we actually want to do is to construct a
group by generators and relations. This problem for geometrical
groups is quite difficult, so we may imagine we already have a
group with such an operation, and this group $G$ is good, it is a
Lie group. Thinking in this way, we see that this group should
have the corresponding Lie algebra $\mathfrak{g}$ (which
\underline{can} in fact be constructed by using generators and
relations), and the group is connected with the algebra by the
exponential mapping $exp: {\mathfrak{g}}\to G$. So, it would be
nice to understand the quandle relations of type (\ref{qrela}),
say, the relation $a\circ b=bab^{-1}$. In the algebra, this
relation would be

\begin{equation}
{\mathfrak{a}}\circ{\mathfrak{b}}=log(exp({\mathfrak{a}})exp({\mathfrak{b}})exp({\mathfrak{c}})),\label{algrel}\end{equation}
where $log$ is the inverse operation $G\to {\mathfrak{g}}$ to the
exponential map (more precisely, this operation is defined in the
vicinity of the group unit, however, we will operate with formulae
and look what happens).

It turns out, that the operation $a,b\mapsto log(exp(a)exp(b))$
can be defined in terms of the Lie algebra, i.e., it can be
expressed in commutators.

The required statement is called the {\em
Baker-Campbell-Hausdorff} theorem, and the coefficients are given
by a beautiful formula due to Dynkin \cite{Dyn}.

\begin{equation} ln(e^{x}e^{y})\sum\frac{(-1)^{k-1}}{k}\frac{1}{{p_1}!{q_1}!\dots
{p_k}!{q_k}!} (x^{p_1}y^{q_1}\dots
x^{p_k}y^{q_k})^{o},\label{BKH}
\end{equation}
where $p_{i}q_{i}>0$ for $i\le k$; here the function
${\cdot}^{\circ}$ is defined on formal non-commutative monomials
in variables $x_{j}$ by the rule

\begin{equation} (x_{i_1}x_{i_2}\dots
x_{i_k})^{\circ}=\frac{1}{k}[\dots[[x_{i_1},x_{i_2}],x_{i_3}],\dots,x_{i_k}]\end{equation}
and extended linearly for their linear combinations.

Thus, for any classical (or virtual) link $L$ we construct an
infinite-dimensional Lie algebra $Li(L)$ as follows (All this can
be done, if we, for instance, take into account virtual crossings
as well; we can just add a new generator ${\mathfrak {q}}$ and
define
$f({\mathfrak{a}})=log(exp({\mathfrak{q}})exp({\mathfrak{a}})exp(-{\mathfrak{q}}))$.
We take its arcs as generators of the free infinite-dimensional
Lie algebra, and factorize the resulting algebra subject to the
relations (\ref{algrel}) taking into account the formula
(\ref{BKH}).)

The obtained infinite-dimensional Lie algebra is a link invariant
is not the most pleasant to work with, however, it admits a
simplification. Thus, for instance, we can take just $a\circ b=
a+[a,b]$.

Decree $a\circ b$ to be $a+[a,b]$ (actually, one should take
$a-[a,b]$, but we choose this operation for simplicity). Then, the
first axiom holds by definition: $[a,a]=a+0=a$. As for the second
axiom, the reverse operation exists by the formula $a\slush
b=a-[a,b]+[[a,b],b]-[[[a,b],b],b]+\dots$ without worrying a lot
whether this formula converges. Even if it has infinitely many
members, we can write all in the finite manner: namely, instead of
writing $p\slush q=r$, we write $r\circ q=p$. There is another way
to handle the situation: we put $a\circ b=a+\eps[a,b]$, where
$\eps$ is just a new generator in the basic ring. In this case
there would be no problem with convergence.

Now, for the third relation we need

$$(a\circ b)\circ c=(a\circ c)\circ (b\circ c).$$

So, we have the relation

$$a+[a,b]+[a,c]+[[a,b],c]=a+[a,c]+[a,b]+[a,[b,c]]+[[a,c],[b,c]]+[[a,b],b]$$
in any of these cases, with or without $\eps$.

Taking into account the Jacobi identity, we see that the equation
above is equivalent to $$[[a,c],[b,c]]=0.$$

This means that we must factorize only by third (not second!)
commutators of the form $[[a,c],[b,c]]$. Thus, we should not cut
our series at some finite place. This leads to some interesting
results.

Denote the obtained invariant (with $a\circ b = a+\eps[a,b]$ at
classical crossings and no structure at virtual crossings) by
$Li(K)$

Moreover, this idea together with the Baker-Campbell-Hausdorff
formula lead to many other realizations and models for quandles.

Having this approach, one may search for particular
representations which certainly lead to new invariants of
classical and virtual knots.

Let us consider the $Li(K)$ invariant (without any structure at
virtual crossings). For instance, the trefoil and the figure eight
knot gives us an finite-dimensional Lie algebra, whence the
(5,2)-torus knot Lie algebra is finite-dimensional.

Indeed, for the trefoil, we have three arcs $a,b,c$ and three
vertices where we have to write down the commutators. Obviously,
we get $a+\eps[a,b]=c$ and two other relations obtained by
permuting this one cyclically. Thus, we can express all first
commutators as linear combinations of generators. So, this algebra
can not be infinite-dimensional.

\begin{figure}
\centering\includegraphics[width=180pt]{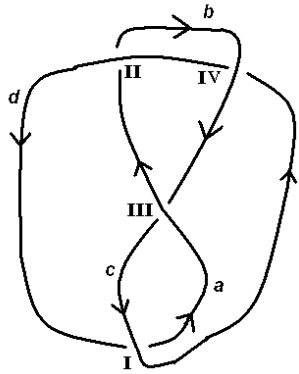}
\caption{The labelled figure eight knot} \label{fig8}
\end{figure}

For the figure eight knot, we have four arcs and three
commutators, one of them is expressed twice. Namely, view Figure
\ref{fig8}.

We see that the first crossing gives us $d\circ c=a$, the second
one gives $b\circ d=a$, the third one gives $b\circ a=c$, and the
fourth one gives $d\circ b=c$. So, we have expressions for three
of six commutators. Thus we have: $\eps[b,d]=a-b,\eps[d,b]=c-d$,
which easily leads to the finite-dimensionality of the algebra.

\begin{figure}
\centering\includegraphics[width=150pt]{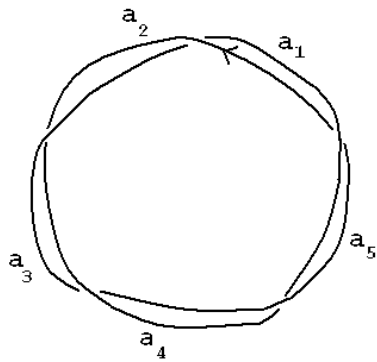} \caption{The torus
(5,2) knot} \label{a5}
\end{figure}

However, if we take the $(5,2)$-torus knot, we have five
generators $a_{1},a_{2},a_{3},a_{4},a_{5}$ and five cyclic
relations $\eps[a_{1},a_{2}]=a_{3}-a_{1},
\eps[a_{2},a_{3}]=a_{4}-a_{2}, \eps[a_{3},a_{4}]=a_{5}-a_{3},
\eps[a_{4},a_{5}]=a_{1}-a_{4}, \eps[a_{5},a_{1}]=a_{2}-a_{5}$, see
Fig. \ref{a5}.

The commutator $[a_{1},a_{3}]$ is not expressible in the terms of
linear combinations of $a_{i}$'s. It is also easy to show that the
elements
$a_{1},([a_{1},a_{3}]),[[a_{1},a_{3}],a_{3}],[[[a_{1},a_{3}],a_{3}],a_{3}]$
represent linearly independent elements. Thus, we obtain a
finite-dimensional algebra.

\section{Colored virtual links \\ and their invariants}
\label{sec}

It is not difficult to show that the virtual quandle construction
generalizes for the case of colored links.

We will not write down all the axioms for the multicomponent
virtual quandle, we are just going to represent one linear model
for oriented virtual links.

Here we give a generalization of the work by Manturov initiated in
\cite{Ma7,Ma9}.

Namely, consider an $n$-component link $L$. Let us fix $n$
generators $t_{1},\dots, t_{n}$ and $n$ generators $s_{1},\dots,
s_{n}$. We are going to associate elements of module over ${\bf
Z}[s_{1},\dots, s_{n},t_{1},\dots, t_{n}]$ to the arcs of the
diagram, whence the module itself will be an invariant. The
relations for virtual crossings look as follows: while passing
through $i$-th component from the left to the right, we multiply
the element associated to the arc ``before'' by $s_{i}$. If we go
from the right to the left, we multiply by corresponding
$s_{j}^{-1}$.

\begin{figure}
\centering\includegraphics[width=180pt]{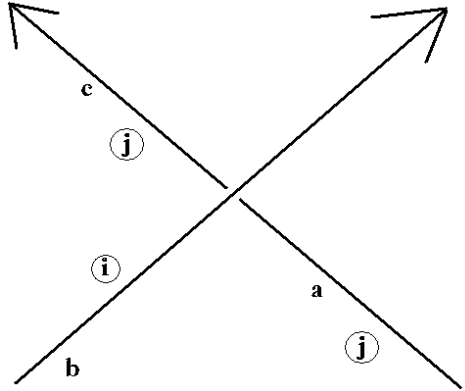} \caption{Relation
at a crossing for colored links} \label{abc}
\end{figure}

At classical crossings, we do the following. Suppose we have a
crossing where $i$-th component goes over, and $j$-th component
goes under, see Fig. \ref{abc}. Suppose the arc going over is $b$,
the arc lying on the right hand is $a$, and that on the left hand
is $c$, see Fig. \ref{abc}.

Then the relation looks like $t_{i}a+(1-t_{j})b=c$.

For negative crossings we write down exactly the same relation
just as in Fig. \ref{abc} with the orientation of $j$-th arc
reversed.

Denote the obtained module by $M(L)$.

The invariance of this module under Reidemeister moves goes
straightforwardly. It turns out that this leads to a polynomial
invariant. This invariant is a generalization of that proposed in
\cite{Ma3} (there we use one variable $t$ and many variables $s$).

\begin{figure}
\centering\includegraphics[width=180pt]{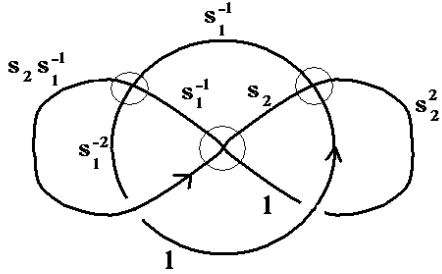}
\caption{Two-component link completely labelled} \label{krug}
\end{figure}

First, we label all arcs of the diagram by monomials in
$s_{1},\dots, s_{n}$ (this will correspond to the operations at
virtual crossings). We do it as follows. First, we deal only with
proper diagrams. Thus, we can associate each crossing with a long
arc. After that, we associate $1$ with each ``first'' arc outgoing
directly from a crossing. After that, we can associate monomials
in $s_{i}$'s to all arcs according to the way described above,
view Fig. \ref{krug}.

Now, we construct a matrix according to the relations at classical
crossings. Namely, each row of our matrix represents a crossing,
each column of the matrix represents a long arc. The matrix
element is going to be the incidence between them. Let us be more
detailed. Suppose we have positive crossing  number $i$ with
incoming edge number $j$ having label $P$ (this edge lies on the
component $p$), and overcrossing edge number $k$ having label $Q$
(this edge lies on the component $q$). Then the $i$-th row of our
matrix should consist of at most three elements, more precisely,
it is equal to the sum of three rows, each of them having only one
non-zero element. One of them has element $i$ equal to $1$,
another one has element $j$ equal to $-t_{q}P$, and the last one
has element $k$ equal to $(t_{p}-1)Q$. In the case when we have
negative crossing, we will have three rows, one of them having
$t_{q}$ on position $i$, another one having $-P$ on position $j$
and the last one having $(t_{p}-1)Q$ on position $k$.

The determinant of the matrix does not change while renumbering
crossings and long arcs correspondingly.

Thus, the determinant is well defined on proper link diagrams. Let
us call this function on diagrams by $\kappa$.

The invariance check is quite similar to that performed in
\cite{Ma3}, so we can show only the most difficult case, namely,
the invariance under the third Reidemeister move.

\begin{figure}
\centering\includegraphics[width=360pt]{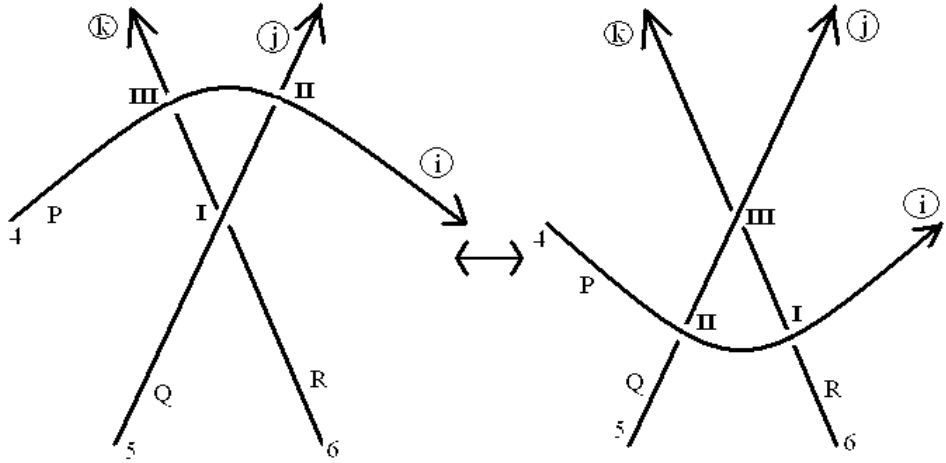}
\caption{Labelling for the third Reidemeister move} \label{lb}
\end{figure}

Namely, consider the two diagram shown in Fig. \ref{lb}.

Here all numbers of crossings are marked by roman numbers, all
numbers of arcs are marked by arabic numbers, numbers of
components are encircled, and monomials corresponding to incoming
arcs are marked by letters $P,Q,R$. View Fig. \ref{lb}.

In the first case we have the matrix

$$\left( \begin{array}{ccccccc}
1 & 0 & 0 & 0 & (t_{k}-1)Q & -t_{j}R & 0\dots 0\\
0 & 1 & 0 & (t_{j}-1)P & -t_{i}Q & 0 & 0\dots 0\\
-t_{i} & 0 & 1 & (t_{k}-1)P & 0 & 0 & 0\dots 0\\
0 \\
\vdots & & & * \\
0
\end{array}\right).$$

The second matrix looks like:

$$\left( \begin{array}{ccccccc}
1 & 0 & 0 & (t_{k}-1)P & 0 & -t_{i}R & 0\dots 0\\
0 & 1 & 0 & (t_{j}-1)P & -t_{i}Q & 0 & 0\dots 0\\
-t_{j} & (t_{k}-1) & 1 & 0 & 0 & 0 & 0\dots 0\\
0 \\
\vdots & & & * \\
0
\end{array}\right).$$

As we see, the second rows of these matrices coincide. So, let us
look at the third rows. For the first matrix,  add the first row
multiplied by $t_{i}$ to the third row. For the second matrix, let
us add the first row multiplied by $t_{j}$ and the second row
multiplied by $(1-t_{k})$, to the third row. We will get two
matrices, for which the third rows coincide. Namely, it looks like

$$(0,0,1,(t_{k}-1)P,t_{i}(t_{k}-1)Q,-t_{i}t_{j}R,0,\dots,0).$$

Now, the first column of the new matrices consists of the only
non-zero element, namely, $1$ on the position $(1,1)$. So, we can
easily make the first rows of these matrices equal. Thus, we have
proved that the determinants of these matrices coincide.

As for the first Reidemeister move, it might multiply the
determinant by a monomial in $t_{i}$'s. Thus, we have proved the
following

\begin{th}
The polynomial $\kappa$ is invariant under Reidemeister moves up
to multiplication by powers of $t_{i}$'s.
\end{th}

The colored link invariant can be used for purposes of knots (not
links) as well: one should just take cabling and care about the
proper normalization of the invariant.

\section{Long Knots and Their Invariants}

It is well known that long classical knots (non-compact knots in
${\bf R}^{3}$ lying on the straight line outside some big circle)
are just the same as ordinary classical knots. Proof see in e.g.
\cite{Ma1}.

The present section the results main ideas of which were sketched
in the book \cite{Ma1}, see also \cite{Ma11}. The two main
arguments that can be taken into account in the theory of the
``long'' virtual knots and could not be used before, are the
following:

\begin{enumerate}

\item One can indicate the {\em initial} and the {\em final} {\em
arcs} (which are not compact) of the diagram representing some two
fixed elements of the quandle; the elements corresponding to them
are invariant under generalized Reidemeister moves.

\item One can take two different quandle--like operations
 at vertices depending on which arc is ``before'' and
which is ``after'' according to the orientation of a long knot.

\end{enumerate}

As shown in \cite{GPV}, the procedure of breaking virtual knot is
not well defined: breaking the same knot diagram at different
points, we may obtain different long knots. Moreover, a virtual
unknot diagram broken at some point can generate a nontrivial long
knot diagram. The aim of this section is to construct invariants
of long virtual knots that ``feel'' the breaking point.

\begin{re}
Throughout the present section, we deal only with long virtual
knots, not links.
\end{re}

\begin{re}
We shall never indicate the orientation of the long knot assuming
it to be oriented from the left to the right.
\end{re}

{\bf Notation}. Throughout this section, $R$ will denote the field
of rational functions in one (real) variable $t$: $R={\bf Q}(t)$.

Now, Let us define the virtual long knot.

\begin{dfn}
By a {\em long virtual knot diagram} we mean a smooth immersion
$f$ of the oriented line $L_{x},x\in(-\infty,+\infty)$ in ${\bf
R}^{2}$, such that:

\begin{enumerate}
\item outside some big circle, we have $f(t)=(t,0)$;

\item each intersection point is double and
transverse;

\item each intersection point is endowed with classical or virtual
crossing structure.
\end{enumerate}
\end{dfn}

\begin{dfn}
A {\em long virtual knot} is an equivalence class of long virtual
knot diagrams modulo generalized Reidemeister moves.
\end{dfn}

\begin{dfn}
A {\em long quandle}\index{Long quandle} is a set $Q$ equipped
with two binary operations $\circ$ and $*$ and one unary operation
$f(\cdot)$ such that $(Q,\circ,f)$ is a virtual quandle
and
$(Q,*,f)$ is a virtual quandle and the following two relations
hold:

The reverse operation for $\circ$ is $\slush$ and the reverse
operation for $*$ is $\slh$.

$$\forall a,b,c\in Q: (a\circ b)*c=(a*c)\circ(b*c),$$
$$\forall a,b,c \in Q: (a*b)\circ c=(a\circ c)*(b\circ c)$$
(new distributivity relations) and
$$\forall x,a,b\in Q: x\alpha (a\circ b)=x\alpha(a*b)$$
$$\forall x,a,b\in Q: x\beta (a\slush b)=x\beta (a\slh b)$$
(strange relations)

where $\alpha$ and $\beta$ are some operations from the list
$\circ,*,\slush,\slh$.
\end{dfn}

\begin{re}
It might seem that the last two relations hold only in
the case when $\circ$ coincides with $*$. However, the
equation $(a\circ b)=c$ has the only relation in $a$ and
not in $b$! As it will be shown later, there are non--trivial
algebraic presentations of the long quandle.
\end{re}

Consider a diagram ${\bar K}$ of a virtual knot and arcs of
it. Let us fix the initial arc $a$ and the
final arc $b$.

Now, we construct the long quandle of it by the following rule.
First, we take all arcs of it including $a$ and $b$ and consider
the {\em free long quandle}: just by using formal operations
$\circ,*,\slush,\slh,f$ factorized only by the quandle relations
(together with the new distributivity relations and strange
relations).

After this, we factorize by relations at crossings. At each
virtual crossing, we do just the same as in the case of virtual
quandle. At each classical crossing we write the relation either
with $\circ$ or with $*$, namely, if the overcrossing is passed
{\bf before} the undercrossing (with respect to the orientation of
the knot) then we use the operation $\circ$ (respectively,
$\slush$); otherwise we use $*$ (respectively, $\slh$).

After this factorization, we obtain an algebraic object $M$
equipped with the five operations $\circ,\slush, *,\slh,$ and $f$
and two selected elements $a$ and~$b$.

\begin{dfn}
Denote the obtained object by $Q_{L}({\bar K})$.
\end{dfn}
Let ${\bar K}$ be a diagram of a long knot $K$. Call $Q_{L}({\bar
K})$ {\em the long quandle} of $K$.

Obviously, for the {\em long unknot} $U$ (represented by a line
without any crossings) we see that the elements $a,b\in Q_{L}(U)$
representing the initial and the final arc should be equal.

\begin{th}
The quandle $Q_{L}$ together with selected elements
$a,b$ is invariant with respect to generalised Reidemeister moves.
\end{th}

\begin{proof}
The proof is quite analogous to the invariance proof of the
virtual quandle.

Thus, the details will be sketched. The invariance under purely
virtual moves and the semivirtual move goes as in the classical
case: we deal only with $f$ and one of the operations $*$ or
$\circ$. Only one of the operators $*,\circ$ appears when applying
the first or the second classical Reidemeister move.

So, the most interesting case is the third classical Reidemeister
move. In fact, it is sufficient to consider the following four
cases shown in Fig. \ref{TRiM} (a,b,c,d). In some cases this means
the distributivity of the operations $\circ,*,\slush,\slh$ with
respect to each other.

\begin{figure}
\centering\includegraphics[width=360pt]{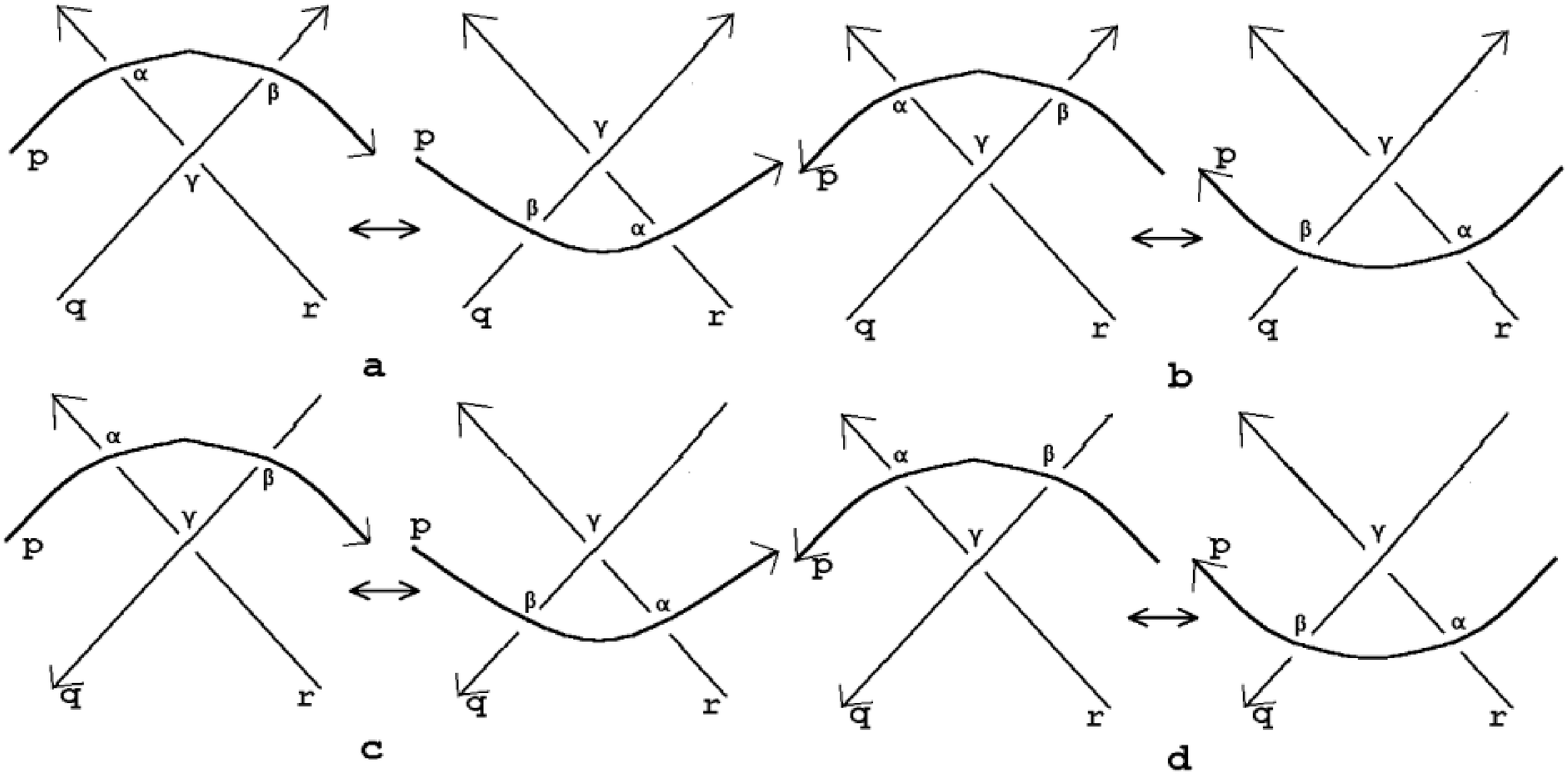}
\caption{Checking the invariance under $\Omega_{3}$} \label{TRiM}
\end{figure}

In each of the four cases everything is OK with $p$ and $q$ ($p$
does not change and $q$ is affected by $p$ in the same manner on
the right hand and on the left hand). So, one should only check
the transformation for $r$.

In each picture, at each crossing we put some operation
$\alpha,\beta$ or $\gamma$. This means one of the operations
$\circ,*,\slush,\slh$ (which are going to be applied to the arc
below to obtain the corresponding arc above).

Consider the case a. We have: each $\alpha,\beta,\gamma$ is a
``multiplications'' $\circ$ or $*$ (the operations $\slush, \slh$
are thus ``divisions'').

Thus, at the upper left corner we shall have: $(r\gamma q)\alpha
p$ in the left picture and $(r\alpha p)\gamma (q\beta p)$ in the
right picture. But, by definition, $(r\gamma q)\alpha p=(r\alpha
p)\gamma (q\alpha p)$. The latter expression equals $(r\alpha
p)\gamma(q\beta p)$ according to the ``new relation'' (because
both $\beta$ and $\alpha$ are ``multiplications'').

Now, let us turn to the case b. Here $\gamma$ is a ``multiplication''
and $\alpha,\beta$ are ``divisions''. Thus, the same equality takes
place: $(r\gamma q)\alpha p=(r\alpha p)\gamma(q\alpha p)=(r\alpha
p)\gamma(q\beta p)$.

The same equality holds for the cases shown in pictures c and d:
the only important thing is that $\alpha$ and $\beta$ are either
both multiplications (as in the case c) or both divisions (as in
the case d). The remaining part of the statement follows
straightforwardly.

\end{proof}

 The (non--trivial) virtual knot represented there is the connected
sum of two unknots. In particular, this means that the
corresponding long virtual knots are not trivial.

\begin{figure}
\centering\includegraphics[width=180pt]{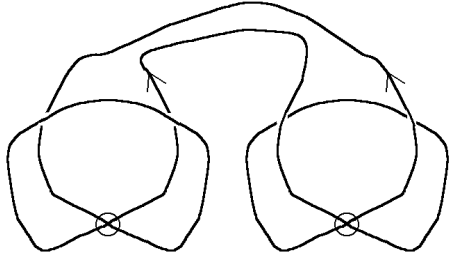} \caption{The
Kishino knot} \label{Kimsg}
\end{figure}

\begin{figure}
\centering\includegraphics[width=320pt]{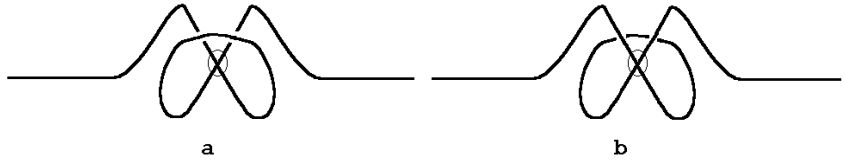} \caption{Two long
virtual knots obtained by breaking the unknot} \label{sgk}
\end{figure}

Consider the unknots shown in Fig. \ref{sgk}, $a$ and $b$. Let us
show that they are not isotopic to the trivial knot. To do it,  we
will use the presentation of the long virtual quandle to the
module over ${\bf Z}_{16}$ by:

$$a\circ b= 5 a-4 b,\; a*b=9a-8b.$$

$$f(x)=3\cdot x.$$

It can be readily checked that these relations satisfy all axioms
of the long quandle.

Let us show that for none of these two knots $a=b$.

Indeed, for the first knot (Fig.\ref{sgk}.a), denote by $c$ the
next arc after $a$. Then we have:

$$9a-8\cdot(3c)=c,5b-4\cdot(3c)=c\Longrightarrow b=9a$$

For the second knot (Fig. \ref{sgk}.b), denote by $c$ the upper
(shortest) arc. We have:

$$5\cdot(3b)-4a=c,9\cdot (3a)-8b=c\Longrightarrow b=9a $$

As we see, in none of these cases $a=b$. Besides, the expressions
of $b$ via $a$ are different. Thus, none of the two long knots
shown in Fig. \ref{sgk}.a and Fig. \ref{sgk}.b is trivial.

If we take another field, say, ${\bf Z}_{25}$ with the operations
$a\circ b=6a-5b$ and $a*b=11a-10b$, $f(x)=3x$ we shall see that
the two knots from Figures \ref{sgk}.a and \ref{sgk}.b are indeed
different: in the first case we obtain $a=11b$, in the second case
we obtain $a=21b$.

In fact, the linear model for long quandles allows to prove even
more, namely, we can show that the knots \ref{sgk}.a and
\ref{sgk}.b do not commute (this implies the non-triviality and
non-classicality of each of them together with their
non-equivalence).

Later, the same effect was established by D.Silver and S.Williams
by using non-commutative structures in biquandles.

\section{Long knot flat quandles}

There is an interesting question to consider: virtual knots modulo
classical crossing change. These objects are called virtual flats,
for more details see, e.g, \cite{Ma1, Ma6}. These object are
classified geometrically, moreover, they lead to powerful
invariants of virtual knots \cite{Ma6,Tur}. However it is worth
studying the algebraic classification and invariants of these
objects (which was performed in \cite{Tur, HK}), because this may
lead to the construction of skein algebras (see e.g. \cite{Ma1})
for virtual knots.

It turns out that this plan can be performed somehow by using the
ideas described in the previous section.

As in \cite{FJK} and \cite{BF}, one can consider linear biquandles
over non-commutative rings, e.g., over quaternions. Namely, having
a classical crossing with two inputs, one writes down two outputs
depending on the first ones linearly by means of some matrix $A$.
Here we do not pay attention to virtual crossings\footnote{This is
an interesting object to be discussed.} Obviously, the second
(classical) Reidemeister move requires invertibility of $A$, and
the third Reidemeister move requires some equation of Yang-Baxter
type, namely,

\begin{equation}A_{1}A_{2}A_{1}=A_{2}A_{1}A_{2},\label{eq8}\end{equation}
where $A_{1}$ is $3\times 3$-matrix consisting of blocks $A$ (of
size $2\times 2$) and $1$ (of size $1\times 1$), and $A_{1}$ is
the matrix consisting of blocks $1$ and $A$.

What have the equations got to do with {\bf long virtual knots}?
It turns out, that we can use two different matrices, say, $A$ and
$B$ for different types of crossings, namely, $A$ for the early
overcrossing and $B$ for the late overcrossing. In this case, we
get the equations (\ref{eq8}) for $A_{1},A_{2}$ and analogous ones
for $B_{1}, B_{2}$ (as before). Besides this, get one more
equation for the third Reidemeister move which corresponds to the
``strange relation'' in the case of long quandles. Namely,

\begin{equation}
A_{1}A_{2}B_{1}=B_{2}A_{1}A_{2}\label{eq9}\end{equation}

Besides obvious solution $B=A$, one can also take the solution
$B=A^{-1}$. It is clear that all condition described above hold
(this is left to the reader as a simple exercise). What do we get
in this case? Obviously, each ``bad crossing'' (with late
overcrossing) is operated on by $B$ (or $B^{-1}$), so, we will
have the same result as in the case of the inverse crossing
($A^{-1}$ or $A$).

Finally, we get the invariant of the ``descending'' long knot with
the same shadow, obtained by using simply the Fenn approach with
the matrix $A$. This is going to be an invariant under generalized
Reidemeister moves, so this is an invariant of flats.

This agrees with the following general statement due to
V.G.Turaev: the mapping associating to a flat long virtual link
diagram the corresponding {\em ascending} long virtual link
diagram is well define. Thus it allows to use any long virtual
knot invariants for recognizing flat long virtual knots.

So, we can take each of the solutions presented in \cite{BF}, and
derive long virtual flat invariants from them.

The question whether these flat long virtual knot invariants can
be used as a basis for some skein module, is still to be
discovered.

\end{document}